\documentclass[11pt,reqno]{amsart}

\usepackage{amsthm}
\usepackage{amssymb}
\usepackage{latexsym}
\usepackage{multicol}
\usepackage{verbatim,enumerate}
\usepackage{accents}
\usepackage[usenames]{color}
\usepackage{youngtab}
\usepackage[all, poly]{xy}  
\usepackage[utf8]{inputenc}
\usepackage{hyperref}
\usepackage{amsmath, amscd}
\usepackage{soul}
\usepackage{tikz} 
\usetikzlibrary{matrix,arrows,decorations.pathmorphing}

\theoremstyle{plain}
\newtheorem*{prop}{Proposition}
\newtheorem{thm}{Theorem}
\newtheorem*{lem}{Lemma}
\newtheorem*{cor}{Corollary}

\theoremstyle{definition}
\newtheorem*{example}{Example}
\newtheorem*{defn}{Definition}

\newtheorem*{rem}{Remark}

\theoremstyle{remark}

\newcommand{\lie}[1]{\mathfrak{#1}}

\newcommand\bz{\mathbb Z}

\newcommand{\qbinom}[2]{\genfrac[]{0pt}0{#1}{#2}}

\advance\textwidth by 1.2in  \advance\oddsidemargin by -.6in \advance\evensidemargin by -.6in

\newcounter{cnt}
 \makeatletter
\def\mydggeometry{\makeatletter\dg@YGRID=1\dg@XGRID=20\unitlength=0.003pt\makeatother}
\makeatother \theoremstyle{remark}

\numberwithin{equation}{section}
\makeatletter
\def\section{\def\@secnumfont{\mdseries}\@startsection{section}{1}%
  \z@{.7\linespacing\@plus\linespacing}{.5\linespacing}%
  {\normalfont\scshape\centering}}
\def\subsection{\def\@secnumfont{\bfseries}\@startsection{subsection}{2}%
  {\parindent}{.5\linespacing\@plus.7\linespacing}{-.5em}%
  {\normalfont\bfseries}}
\makeatother

\begin{document}

\title[Prime representations in the HL category: classical decompositions]{Prime representations in the Hernandez-Leclerc category: classical decompositions}
\author{Leon Barth}
\address{Faculty of Mathematics Ruhr-University Bochum}
\email{leon.barth@rub.de}
\author{Deniz Kus}
\address{Faculty of Mathematics Ruhr-University Bochum}
\email{deniz.kus@rub.de}
\thanks{D.K. was partially funded by the Deutsche Forschungsgemeinschaft (DFG, German Research Foundation) – 446246717.}

\subjclass[2010]{}
\begin{abstract} We use the dual functional realization of loop algebras to study the prime irreducible objects in the Hernandez-Leclerc category for the quantum affine algebra associated to $\mathfrak{sl}_{n+1}$. When the HL category is realized as a monoidal categorification of a cluster algebra \cite{HL10,HL13}, these representations correspond precisely to the cluster variables and the frozen variables are minimal affinizations. For any height function, we determine the classical decomposition of these representations with respect to the Hopf subalgebra $\mathbf{U}_q(\mathfrak{sl}_{n+1})$ and describe the graded multiplicities of their graded limits in terms of lattice points of convex polytopes. Combined with \cite{BCMo15} we obtain the graded decomposition of stable prime Demazure modules in level two integrable highest weight representations of the corresponding affine Lie algebra. 

\end{abstract}

\maketitle

\section{Introduction}
The classification of finite-dimensional irreducible representations of quantum affine algebras was given in \cite{CP91,CP95} in terms of Drinfeld polynomials. However, generically the structure of these representations is far from beeing understood. Not even dimension formulas exist in contrast to the classical cases except for particular families of modules, e.g. local Weyl modules, Kirillov–Reshetikhin modules or minimal affinizations (see \cite{CP01,He06,He07,Nak03} for instance).
A well established method to study these representations is to pass from quantum level to classical level by forming their classical limit; see for instance \cite{CP01} for a necessary and sufficient condition for the existence of the classical limit. When the limit exists, it is a finite–dimensional module for the corresponding affine Lie algebra and hence a representation for the underlying standard maximal parabolic subalgebra - the current algebra. After a suitable twist, which is referred to as the graded limit in the literature, many interesting families of graded representations of current algebras appear in this way. Among them are the Kirillov-Reshetikhin modules for current algebras (see \cite{CM04,FK08,Ke11}) and their fusion products \cite{N17a} and also generalized Demazure modules appear in this context \cite{CDM19}. Inspired by the results of \cite{AK07,AKS06,N17a} the motivation of this paper is to fully understand the graded decompositions of the graded limits of the prime irreducible objects in the Hernandez-Leclerc category which we explain now in more detail.

Let $\mathfrak{g}$ be the finite-dimensional complex simple Lie algebra $\mathfrak{sl}_{n+1}$ of $(n+1)\times (n+1)$ complex matrices of trace zero. In the seminal paper of Hernandez and Leclerc \cite{HL10} the authors presented an interesting subcategory $\mathcal{C}_{q,\kappa}$ of the category of all finite-dimensional representations depending on a height function $\kappa$. Their main result states that $\mathcal{C}_{q,\kappa}$ is closed under tensor products and categorifies a cluster algebra of the same type. The prime irreducible objects in that category, i.e. the ones which are not isomorphic to a tensor product of non–trivial representations, correspond  precisely to the cluster variables and the frozen variables correspond to minimal affinizations. Moreover, an explicit description of the prime objects in terms of Drinfeld polynomials is given in \cite{BC19a,HL13} (see also Theorem~\ref{hlmaan}) and their graded limits are isomorphic to stable prime Demazure modules in level two integrable highest weight representations \cite{BCMo15}. In this paper we give a description of the structure of these objects viewed as representations for the Hopf subalgebra $\mathbf{U}_q(\mathfrak{sl}_{n+1})$. We explain now the results of the paper in more detail. 

For a finite-dimensional graded representation $V$ for the current algebra we denote by $\tau_p^{*}V$ the graded vector space whose $r$-th graded piece is $V[r+p]$. We encode the graded multiplicities as polynomials in an indeterminate $q$ as follows
$$\big[V: V(\mu)\big]_q= \sum_{p=0}^{\infty}\ [V: \tau_p^{*}V(\mu)] \cdot q^p.$$
Our aim is to determine these polynomials for a wide class of graded  representations. In Section~\ref{section4} we introduce the modules $M_{\boldsymbol \xi,\lambda}$ depending on a pair $(\boldsymbol \xi,\lambda)$ where $\lambda$ is a dominant integral weight and $\boldsymbol \xi$ is a tuple of non-negative integers indexed by the positive roots of $\mathfrak{g}$. Special choices of $M_{\boldsymbol \xi,\lambda}$ give many well-known families of representations such as truncated Weyl modules (see \cite{BK20,FMM19,KL14} for instance) or the graded limits of certain representations for quantum affine algebras. We first describe the $\mathbf{U}^{-}$ structure of these representations in Proposition~\ref{relc} by generators and relations. Subsequently, using the dual functional realization of loop algebras and the methods developed in \cite{AK07,AKS06} we give a functional description of the graded multiplicities of $M_{\boldsymbol \xi,\lambda}$ in Theorem~\ref{wicht}. This space of functions can be always identified with a subalgebra (sometimes a representation) of the Cohomological Hall algebra (CoHA) of a quiver which was introduced by Kontsevich and Soibelman in \cite{KS11a}. As of now, there are very few examples of modules or of subalgebras of the CoHA, so that this connection could be of independent interest; the details will appear elsewhere.

Our goal is to figure out more explicit descriptions of these multiplicities, for example combinatorial parametrizations, however this question seems to be quite challenging for arbitrary pairs $(\boldsymbol \xi,\lambda)$. One possible explanation is the following. If $\boldsymbol{\xi}$ is a constant tuple, say each entry is equal to $N$, and $\lambda=N\mu$ for a dominant integral weight $\mu$, then the numerical multiplicity $[M_{\boldsymbol{\xi},\lambda}: V(\nu)]_{q=1}$ is exactly the Littlewood-Richardson coefficient describing how often $V(\nu)$ appears inside the tensor product $V(\mu)^{\otimes N}$ (see \cite{KL14} for instance). 

The main result of the paper gives a description of the graded multiplicities (see Theorem~\ref{mainthm}) of the graded limits $L(\boldsymbol{\pi})$ in terms of lattice points of convex polytopes when $\boldsymbol{\pi}$ is the Drinfeld polynomial whose corresponding representation $V_q(\boldsymbol{\pi})$ is a prime irreducible object in $\mathcal{C}_{q,\kappa}$. In fact the result can be applied for all representations whose generators and relations have a particular form (see Remark~\ref{partform}), e.g. for certain generalized Demazure modules or minimal affinizations by parts.
 
\textit{Organization of the paper:} In Section~\ref{section2} we introduce the main definitions and notations and discuss the prime irreducible objects in the Hernandez-Leclerc category as well as their graded limits. In Section~\ref{section3} we state the main results and in Section~\ref{section4} we present a class of truncated representations for current algebras and determine their $\mathbf{U}^-$ module structure. In Section~\ref{section5} we recall the dual functional realization of loop algebras and discuss the graded characters of the aforementioned truncated representations. In Section~\ref{section6} we prove the main theorem of the paper giving graded decompositions of the prime irreducible objects in the HL category in terms of lattice points of convex polytopes.

\textit{Acknowledgement: D.K. thanks the Hausdorff Research Institute for Mathematics and the organizers of the Trimester Program ``New Trends in Representation Theory"  for excellent working conditions. He also thanks Matheus Brito for many useful discussions and Rekha Biswal for drawing our attention to the dual functional realzation which was one of the key steps in proving the results and also for many inspiring discussions especially in Sections~\ref{section4}-\ref{section5}.}
\section{Quantum loop algebras and prime representations in the HL category}\label{section2}
\subsection{}
Throughout this paper we denote by $\mathbb{C}$ the field of complex numbers and by $\mathbb{Z}$ (resp. $\mathbb{Z}_{+}$, $\mathbb{N}$) the subset of integers (resp. non-negative, positive integers).
\subsection{} Let $\mathfrak{g}$ the Lie algebra $\mathfrak{sl}_{n+1}$ of $(n+1)\times (n+1)$ complex matrices of trace zero with Cartan matrix $(c_{i,j})$. Let $\mathfrak{h}$ be a Cartan subalgebra, $R$ be the corresponding set of roots and $\{\alpha_1,\dots,\alpha_n\}$ and $\{\varpi_1,\dots,\varpi_n\}$ respectively be a set of simple roots and fundamental weights respectively. We denote the set of positive roots by $R^+$, the $\mathbb{Z}$ (resp. $\mathbb{Z}_+$) span of the simple roots by $Q$ (resp. $Q^+$) and the $\mathbb{Z}$ (resp. $\mathbb{Z}_+$) span of the fundamental weights by $P$ (resp. $P^+$). We define as usual a partial order on $P$ by $\lambda \geq \mu $ if $\lambda-\mu\in Q^+$. Note that for two positive roots $\alpha,\beta\in R^+$ with $\beta> \alpha$ we either have $\beta-\alpha\in R^+$ or there exists $\gamma_1,\gamma_2\in R^+$ such that 
\begin{equation}\label{seqr}\beta-\gamma_1\in R^+,\ \ \alpha=\beta-\gamma_1-\gamma_2.\end{equation}
For a given root $\alpha\in R^+$, let $x_{\alpha}^{\pm}$ be the corresponding root vector of weight $\pm \alpha$ and $h_{\alpha}$ the corresponding coroot. We have a triangular decomposition 
$$\mathfrak{g}=\mathfrak{n}^{-}\oplus\mathfrak{h}\oplus \mathfrak{n}^{+},\ \ \mathfrak{n}^{\pm}=\bigoplus_{\alpha\in R^+} \mathbb{C} \cdot x_{\alpha}^{\pm}.$$
 If $\alpha=\alpha_i+\cdots+\alpha_j$, $1\leq i\leq j\leq n$, and $\lambda\in P^+$ we abbreviate in the rest of the paper  $$x^{\pm}_{\alpha}:=x^{\pm}_{i,j},\ x^{\pm}_{i,i}:=x^{\pm}_i,\ h_{\alpha}:=h_{i,j},\ h_{i,i}:=h_i,\ \lambda(h_{\alpha}):=\lambda_{i,j},\ \lambda_{i,i}:=\lambda_i.$$
The unique irreducible representation of $\mathfrak{g}$ of highest weight $\lambda\in P^+$ is denoted by $V(\lambda)$.
\subsection{}Let $\widehat{\mathfrak{g}}$ the untwisted affine Lie algebra associated to $\mathfrak{g}$ which is realized as 
$$\widehat{\mathfrak{g}}=\mathfrak{g}\otimes\mathbb{C}[t^{\pm}]\oplus \mathbb{C}K\oplus \mathbb{C}d$$ where $K$ is required to be central and the Lie bracket is defined as
$$[x\otimes t^{r},y\otimes t^s]=[x,y]\otimes t^{r+s}+\mathrm{tr}(xy)K,\ \ [d,x\otimes t^r]=r(x\otimes t^r),\ \ \ x,y\in\mathfrak{g},\ \ r,s\in\mathbb{Z}.$$
The commutator subalgebra $[\widehat{\mathfrak{g}},\widehat{\mathfrak{g}}]$ modulo the center is the loop algebra $\mathfrak{g}\otimes \mathbb{C}[t^{\pm}]:=\mathfrak{g}[t^{\pm}]$ and note that the element $d$ defines a grading on the loop algebra. The $\mathbb{Z}_+$-graded subalgebra $\mathfrak{g}\otimes \mathbb{C}[t]$ of the loop algebra is the current algebra associated to $\mathfrak{g} $ which we denote by $\mathfrak{g}[t]$. 
For a Lie algebra $\mathfrak{a}$, let $\mathbf{U}(\mathfrak{a})$ be the universal enveloping algebra of $\mathfrak{a}$. In the rest of the paper we abbreviate  
$$\mathbf{U}=\mathbf{U}(\mathfrak{g}[t]),\ \ \mathbf{U}^{\pm}=\mathbf{U}(\mathfrak{n}^{\pm}[t]),\ \ \mathbf{U}^0=\mathbf{U}(\mathfrak{h}[t]).$$ 

$$\mathbf{U}_{\ell}=\mathbf{U}(\mathfrak{g}[t^{\pm}]),\ \ \mathbf{U}_{\ell}^{\pm}=\mathbf{U}(\lie n^{\pm}[t,t^{-1}]),\ \ \mathbf{U}^0_{\ell}=\mathbf{U}(\mathfrak{h}[t^{\pm}])$$
So as vector spaces,
$$\mathbf{U}\cong \mathbf{U}^-\otimes \mathbf{U}^{0} \otimes \mathbf{U}^+,\ \ \mathbf{U}_{\ell}\cong \mathbf{U}_{\ell}^-\otimes \mathbf{U}_{\ell}^{0} \otimes \mathbf{U}_{\ell}^+.$$
\subsection{}Let $\mathbb{C}(q)$ be the field of rational functions in an indeterminate $q$. We discuss in the rest of this section quantum loop algebras, their representations (of type 1) and graded limits. Set
$$[m]=\frac{q^{m}-q^{-m}}{q-q^{-1}},\ \ [m]!=[m][m-1]\cdots [1],\ \ \qbinom{m}{r}=\frac{[m]!}{[r]![m-r]!},\ \ \ r,m\in \mathbb{Z}_+,\ m\geq r.$$
The quantum loop algebra $\mathbf{U}_q(\mathfrak{g}[t^{\pm}])$ is the $\mathbb{C}(q)$-algebra generated by elements $$\tilde{x}_{i,r}^{{}\pm{}},\ \tilde{k}_i^{{}\pm 1},\ \tilde{h}_{i,s}\ \ (1\leq i\leq n, r\in \mathbb{Z}, s\in \mathbb{Z}\backslash\{0\})$$ subject to the following relations:
$$\tilde{k}_i\tilde{k}_i^{-1} = \tilde{k}_i^{-1}\tilde{k}_i =1,\ \ [\tilde{k}_i,\tilde{k}_j]=[\tilde{k}_i,\tilde{h}_{j,r}] =[\tilde{h}_{i,r},\tilde{h}_{j,s}]=0,\ \
\tilde{k}_i\tilde{x}_{j,r}^{{}\pm{}}\tilde{k}_i^{-1} = q^{{}\pm c_{ij}}\tilde{x}_{j,r}^{{}\pm{}},$$
$$[\tilde{h}_{i,r} , \tilde{x}_{j,s}^{{}\pm{}}] = \pm\frac1r[rc_{ij}]\tilde{x}_{j,r+s}^{{}\pm{}},$$
$$\tilde{x}_{i,r+1}^{{}\pm{}}\tilde{x}_{j,s}^{{}\pm{}} -q^{{}\pm
a_{ij}}\tilde{x}_{j,s}^{{}\pm{}}\tilde{x}_{i,r+1}^{{}\pm{}} =q^{{}\pm
a_{ij}}\tilde{x}_{i,r}^{{}\pm{}}\tilde{x}_{j,s+1}^{{}\pm{}}
-\tilde{x}_{j,s+1}^{{}\pm{}}\tilde{x}_{i,r}^{{}\pm{}},$$
$$[\tilde{x}_{i,r}^{{}\pm{}},\tilde{x}_{j,s}^{{}\pm{}}]=0,\ \text{if $c_{i,j}=0$}$$
$$\tilde{x}_{j,s}^{{}\pm{}}
 \tilde{x}_{i, r_{1}}^{{}\pm{}}\tilde{x}_{i,r_{2}}^{{}\pm{}}+\tilde{x}_{j,s}^{{}\pm{}}
 \tilde{x}_{i, r_{2}}^{{}\pm{}}\tilde{x}_{i,r_{1}}^{{}\pm{}}+x_{i, r_{1}}^{{}\pm{}}\tilde{x}_{i,r_{2}}^{{}\pm{}}\tilde{x}_{j,s}^{{}\pm{}}+x_{i, r_{2}}^{{}\pm{}}\tilde{x}_{i,r_{1}}^{{}\pm{}}\tilde{x}_{j,s}^{{}\pm{}}=[2](\tilde{x}_{i, r_{1}}^{\pm}
\tilde{x}_{j,s}^{{}\pm{}}\tilde{x}_{i, r_{2}}^{\pm}+\tilde{x}_{i, r_{2}}^{\pm}
\tilde{x}_{j,s}^{{}\pm{}}
 \tilde{x}_{i, r_{1}}^{{}\pm{}}),\ i\neq j$$
 $$[\tilde{x}_{i,r}^+ , \tilde{x}_{j,s}^-]=\delta_{ij}  \frac{\phi_{i,r+s}^+ -
\phi_{i,r+s}^-}{q - q^{-1}},$$
where  $\phi_{i,r}^{\pm}$ is determined by equating coefficients of powers of $u$ in 
$$\Phi_i^{\pm}(u)=\sum_{r\in\mathbb{Z}}\phi_{i,\pm r}^{{}\pm{}}u^{{} r} = \tilde{k}_i^{{}\pm 1}
\mathrm{exp}\left(\pm(q-q^{-1})\sum_{s=1}^{\infty}\tilde{h}_{i,\pm s} u^{{}s}\right).$$
Denote by $\mathbf{U}_q(\lie g)$ and $\mathbf{U}_q(\lie h[t^{\pm}])$ respectively the subalgebra generated by $\{\tilde{x}_{i,0}^\pm, \tilde{k}_i^{\pm 1}, 1\leq i\leq n\}$ and $\{\tilde{k}_i^{\pm1},
\tilde{h}_{i,s}, 1\leq i\leq n, s\in\mathbb{Z}\backslash\{0\}\}$ respectively. Then $\Lambda_{i,r}$ together with $\tilde{k}_i^{\pm1}$, $1\leq i\leq n$, $r\in \mathbb{Z}$, generate $\mathbf{U}_q(\lie h[t^{\pm}])$ as an algebra where the elements $\Lambda_{i,r}$ are obtained by equating powers of $u$ in the formal power series
\begin{equation*}
\Lambda_i^\pm(u)=\sum_{r=0}^\infty \Lambda_{i,\pm r} u^{r}=
\exp\left(-\sum_{s=1}^\infty\frac{\tilde{h}_{i,\pm s}}{[s]}u^s\right).
\end{equation*}
\subsection{} We consider the dominant $\ell$-weight lattice of $\mathbf{U}_q(\mathfrak{g}[t^{\pm}])$ defined as the monoid $\mathcal P^+$ of $n$-tuples of polynomials $\boldsymbol{\pi} = (\boldsymbol{\pi}_1(u),\dots,\boldsymbol{\pi}_n(u))$ with coefficients in $\mathbb{C}(q)[u]$ such that $\boldsymbol{\pi}_i(0)=1$ for all $1\leq i\leq n$. Given $a\in\mathbb{C}(q)^{\times}$ and $1\leq i\leq n$, define the fundamental $\ell$-weight $\boldsymbol{\varpi}_{i,a}\in \mathcal P^+$ by
$$(\boldsymbol{\varpi}_{i,a})_j(u) = (1-\delta_{i,j}au).$$
The $\ell$-weight lattice $\mathcal P$ is the free abelian group generated by fundamental $\ell$-weights. We consider the map $\boldsymbol{\Psi}:\mathcal{P}\rightarrow  \mathbf{U}_q(\lie h[t^{\pm}])^{*},\ \boldsymbol{\pi}\rightarrow \boldsymbol{\Psi}_{\boldsymbol{\pi}}$ (which turns out to be injective) by the following rule on the generators. For $\boldsymbol{\pi}=\boldsymbol{\pi}'\widetilde{\boldsymbol{\pi}}^{-1}$ with $\boldsymbol{\pi}',\widetilde{\boldsymbol{\pi}}\in\mathcal{P}^+$, let $\boldsymbol{\Psi}_{\boldsymbol{\pi}}$ the algebra homomorphism determined by
$$\boldsymbol{\Psi}_{\boldsymbol{\pi}}(\tilde{k}_i^{\pm1})=q^{\pm \mathrm{wt}(\boldsymbol{\pi})(h_i)},\ \ \boldsymbol{\Psi}_{\boldsymbol{\pi}}(\Lambda^{\pm}_{i}(u))=\frac{\boldsymbol{\pi}'^{\pm}_i(u)}{\widetilde{\boldsymbol{\pi}}^{\pm}_i(u)}$$
where the weight map ist the group homomorphism $\mathrm{wt}:\mathcal P \to P$,\ $\mathrm{wt}(\boldsymbol\varpi_{i,a})=\varpi_i$ and $\boldsymbol{\pi}^{+}_i(u)=\boldsymbol{\pi}^{}_i(u)$ whereas $\boldsymbol{\pi}^{-}_i(u)$ is the polynomial obtained from $\boldsymbol{\pi}_i(u)$ by replacing each $\boldsymbol{\varpi}_{i,a}$ by $\boldsymbol{\varpi}_{i,a^{-1}}$.
\subsection{} A nonzero vector $v$ of a $\mathbf{U}_q(\lie g[t^{\pm}])$-module $V$ 
is called an $\ell$-weight vector of $\ell$-weight $\boldsymbol{\pi}\in\mathcal P$ if there exists $k\in\mathbb{N}$ such that
$$(H-\boldsymbol{\Psi}_{\boldsymbol{\pi}}(H))^kv=0 \ \ \text{ for all }\ \ H\in \mathbf{U}_q(\lie h[t^{\pm}]).$$ If we have 
$$Hv=\boldsymbol{\Psi}_{\boldsymbol{\pi}}(H)v \ \ \text{for all}\ \ H\in \mathbf{U}_q(\lie h[t^{\pm}])
\ \ \text{and}\ \ \tilde{x}_{i,r}^+v=0 \quad\text{for all}\quad 1\leq i\leq n,\ r\in\mathbb Z,$$
then $v$ is called a highest $\ell$-weight vector. The representation $V$ is called an $\ell$-weight module if every vector of
$V$ is a linear combination of $\ell$-weight vectors and a highest $\ell$-weight module if it is generated by a
highest $\ell$-weight vector. Let $\mathcal C_q$ be the category of
all finite-dimensional $\ell$-weight modules of $\mathbf{U}_q(\lie g[t^{\pm}])$.
Note that $\mathcal C_q$ is an abelian category stable under tensor product; see for example \cite{FR99}. The irreducible objects were classified by Chari and Pressley and are obtained as follows. Let $\boldsymbol{\pi}\in\mathcal{P}^+$ and $W_q(\boldsymbol{\pi})$ be the $\mathbf{U}_q(\lie g[t^{\pm}])$-module generated by an element $w_{\boldsymbol{\pi}}$ with defining relations
$$\tilde{x}_{i,r}^+w_{\boldsymbol{\pi}}=(x_{i,0}^-)^{\mathrm{wt}(\boldsymbol{\pi})(h_i)+1}w_{\boldsymbol{\pi}}=(H-\boldsymbol\Psi_{\boldsymbol{\pi}}(H))w_{\boldsymbol{\pi}}=0,\  1\leq i\leq n,\ r\in\mathbb Z,\ H\in \mathbf{U}_q(\lie h[t^{\pm}]).$$
Since $W_q(\boldsymbol{\pi})$ is a highest $\ell$-weight module, it has a unique irreducible quotient $V_q(\boldsymbol{\pi})$. The next theorem can be derived from \cite{CP91,CP95}.
\begin{thm} The map $$\mathcal{P}^+\rightarrow \{\text{irreducible objects in $\mathcal{C}_q$}\}/\sim$$$$ \boldsymbol{\pi}\mapsto V_q(\boldsymbol{\pi})$$
is a one-to-one correspondence between the irreducible objects in the category $\mathcal{C}_q$ and the dominant $\ell$-weights.
\hfill\qed
\end{thm}
\subsection{}There are several ways to study these representations. One method would be to determine the classical decompostion with respect to the Hopf subalgebra $\mathbf{U}_q(\mathfrak{g})$. If $V\in\mathcal{C}_q$, then $V$ can be viewed as a $\mathbf{U}_q(\mathfrak{g})$-module and hence (the category of finite-dimensional $\mathbf{U}_q(\mathfrak{g})$-modules is semisimple)
$$V\cong \bigoplus_{\mu\in P^+} V_q(\mu)^{c_{\mu}}$$
where $V_q(\mu)$ is the irreducible highest weight module of $\mathbf{U}_q(\mathfrak{g})$ of highest weight $\mu$. To be more precise, $V_q(\mu)$
is generated by a vector $v$ satisfying
$$\tilde{x}_{i,0}^+v=0, \quad \tilde{k}_iv=q^{\mu(h_i)}v, \quad (\tilde{x}^{-}_{i,0})^{\mu(h_i)+1}v=0,\quad\forall\ 1\leq i\leq n.$$
We describe now the graded limit approach to the irreducible objects in $\mathcal{C}_q$. Let $\mathbf{A}=\mathbb{Z}[q^{\pm}]$ and denote by $\mathbf{U}_{q,\mathbf{A}}(\mathfrak{g})$ and $\mathbf{U}_{q,\mathbf{A}}(\mathfrak{g}[t^{\pm}])$ respectively the $\mathbf{A}$-form of $\mathbf{U}_{q}(\mathfrak{g})$ and $\mathbf{U}_{q}(\mathfrak{g}[t^{\pm}])$ respectively (for a precise definition see \cite{L2010}). These are free modules over the ring $\mathbf{A}$ such that 
$$\mathbf{U}_{q}(\mathfrak{g})\cong \mathbf{U}_{q,\mathbf{A}}(\mathfrak{g})\otimes_{\mathbf{A}} \mathbb{C}(q),\ \ \mathbf{U}_{q}(\mathfrak{g}
[t^{\pm}])\cong \mathbf{U}_{q,\mathbf{A}}(\mathfrak{g}
[t^{\pm}])\otimes_{\mathbf{A}} \mathbb{C}(q).$$ 
We could try to mimic the same kind of construction for an arbitrary irreducible object in $\mathcal{P}^+$, but the existence is not guaranteed in general. Assume that $V_q(\boldsymbol{\pi})$ admits an $\mathbf{A}$-form, i.e. there is a representation $V_{q,\mathbf{A}}(\boldsymbol{\pi})$ of $\mathbf{U}_{q,\mathbf{A}}(\mathfrak{g}[t^{\pm}])$ such that 
$$V_{q}(\boldsymbol{\pi})\cong V_{q,\mathbf{A}}(\boldsymbol{\pi})\otimes_{\mathbf{A}} \mathbb{C}(q).$$
The classical limit is defined as 
$$\overline{V_{q}(\boldsymbol{\pi})}:= V_{q,\mathbf{A}}(\boldsymbol{\pi})\otimes_{\mathbf{A}} \mathbb{C},$$
where $\mathbb{C}$ is viewed as an $\mathbf{A}$-module by letting $q$ act as $1$.
Since $\mathbf{U}_{q,\mathbf{A}}(\mathfrak{g}[t^{\pm}])\otimes_{\mathbf{A}} \mathbb{C}$ is a quotient of the universal enveloping algebra $\mathbf{U}_{\ell}$, we obtain that $\overline{V_{q}(\boldsymbol{\pi})}$ is a module for $\mathbf{U}_{\ell}$ and hence for $\mathbf{U}$ by restriction. The graded limit $L(\boldsymbol{\pi})$ is obtained by pulling back the $\mathbf{U}$-module $\overline{V_{q}(\boldsymbol{\pi})}$ via the automorphism 
$$\mathfrak{g}[t]\rightarrow \mathfrak{g}[t],\ \ x\otimes t^r\mapsto x\otimes (t-1)^r.$$
So whenever $V_q(\boldsymbol{\pi})$ admits an $\mathbf{A}$-form, we can associate a representation $L(\boldsymbol{\pi})$ of $\mathbf{U}$ to it.
\subsection{} We follow the notation of \cite{BCMo15} and introduce the subset $\mathcal{P}^+_{\mathbb{Z}}(1)$. Let $\mathcal{P}^+_{\mathbb{Z}}$ be the submonoid of $\mathcal{P}^+$ generated by the elements $\boldsymbol{\varpi}_{i,a}$ with $a\in q^{\mathbb{Z}}$. Then it has been shown (see for example the results of \cite{CP01}) that $V_q(\boldsymbol{\pi})$ admits an $\mathbf{A}$-form for all $\boldsymbol{\pi}\in \mathcal{P}^+_{\mathbb{Z}}$ and hence the graded limit exists. Furthermore, we denote
by $\mathcal{P}_{\mathbb{Z}}^+(1)$ the subset of $\mathcal{P}^+_{\mathbb{Z}}$ consisting of elements 
$$\boldsymbol{\varpi}_{i_1,a_1}\cdots \boldsymbol{\varpi}_{i_k,a_k},\ \ 1\leq i_1<i_2\cdots<i_k\leq n,\ \ a_j\in q^\mathbb{Z}$$
such that 
$$a_{j+1}=a_jq^{\pm(i_{j+1}-i_j+2)}$$
and 
$$a_{j+1}=a_jq^{\pm(i_{j+1}-i_j+2)}\Rightarrow a_{j+2}=a_{j+1}q^{\mp(i_{j+2}-i_{j+1}+2)}.$$
The graded limit $L(\boldsymbol{\pi})$ has been determined in \cite{BCMo15} and is isomorphic to a level two Demazure module of weight $\mathrm{wt}(\boldsymbol{\pi})$. To be more precise the following result has been proved.
\begin{thm}\label{mainbcm} Let $\boldsymbol{\pi}=\boldsymbol{\varpi}_{i_1,a_1}\cdots \boldsymbol{\varpi}_{i_k,a_k}\in \mathcal{P}^+_{\mathbb{Z}}(1)$ and $\lambda=\mathrm{wt}(\boldsymbol{\pi})$. Then $L(\boldsymbol{\pi})$ is isomorphic to the graded $\mathbf{U}$-module generated by an element $v$ of grade zero with defining relations:
$$\mathfrak{n}^+[t]v=0,\ (h_i\otimes t^r)v=\delta_{r,0} \lambda_i\ (r\in\mathbb{Z}_+,\ 1\leq i\leq n),$$
$$(x_i^{-}\otimes 1)^{\lambda_i+1}v=0\ (1\leq i\leq n),\ (x^{-}_{i_j,i_{j+1}}\otimes t)v=0,\ 1\leq j\leq k-1.$$ 
\end{thm}
We describe now the importance of these representations and the relation to the Hernandez-Leclerc category. 
\begin{defn}
We call $V_q(\boldsymbol{\pi})$ a prime irreducible representation of $\mathbf{U}_{q}(\mathfrak{g}[t^{\pm}])$ if 
$$V_q(\boldsymbol{\pi})\cong V_q(\boldsymbol{\pi}^1)\otimes \cdots \otimes V_q(\boldsymbol{\pi}^s)$$ implies that $(s-1)$ factors are trivial representations.
\end{defn}
It is clear that $V_q(\boldsymbol{\pi})$ is either prime or can be written as a tensor product of non–trivial prime representations; however the uniqueness of such a decmposition is not known in general. The motivation of this paper is to determine the structure of the prime objects in the HL category by describing the graded character of their graded limits. First we recall the definition. Let $\kappa: \{1,\dots,n\}\rightarrow \mathbb{Z}$ be a height function satisfying $|\kappa(i+1)-\kappa(i)|\leq 1$ for $1\leq i\leq n$ and let $Q_{\kappa}$ the corresponding quiver whose vertices are indexed by $\{1,\dots,n\}$ and there is an edge $i\rightarrow i+1$ if $\kappa(i)<\kappa(i+1)$ and $i\leftarrow i+1$ otherwise. The Hernandez-Leclerc category $\mathcal{C}_{q,\kappa}$ is the full subcategory of $\mathcal{C}_q$ whose objects have all its Jordan-Hölder components of the form 
$$V_q(\boldsymbol{\pi}),\ \ \boldsymbol{\pi}\in \mathcal{P}_{\mathbb{Z}}^+(\kappa,1)$$
where $\mathcal{P}_{\mathbb{Z}}^+(\kappa,1)$ is the submonoid of $\mathcal{P}^+$ generated by $\boldsymbol{\varpi}_{i,a},\ a\in \{q^{\kappa(i)},q^{\kappa(i)+2}\}, 1\leq i\leq n.$ The following results have been proved in \cite{BC19a,HL10, HL13}.
\begin{thm}\label{hlmaan} The category $\mathcal{C}_{q,\kappa}$ is closed under tensor products. Let $V_q(\boldsymbol{\pi})$ be a prime irreducible object in $\mathcal{C}_{q,\kappa}$. Then $\boldsymbol{\pi}\in\{\boldsymbol{\varpi}_{i,q^{\kappa(i)}}\boldsymbol{\varpi}_{i,q^{\kappa(i)+2}},\boldsymbol{\varpi}_{i,q^{\kappa(i)+2}},\boldsymbol{\varpi}_{i,q^{\kappa(i)}}\}$ for some $i\in\{1,\dots,n\}$ or there exists an interval $J\subseteq[1,n]$ such that 
$$\boldsymbol{\pi}=\boldsymbol{\pi}_{\kappa,J}:=\prod_{i\in J_{\mathrm{sink}}} \boldsymbol{\varpi}_{i,q^{\kappa(i)}}\prod_{i\in J_{\mathrm{source}}} \boldsymbol{\varpi}_{i,q^{\kappa(i)+2}}$$
where $J_{\mathrm{sink}}$ (resp. $J_{\mathrm{source}}$) are the sinks (resp. sources) of $Q_{\kappa}$ contained in $J$. Conversely, all these representations are prime objects in $\mathcal{C}_{q,\kappa}$.
\hfill\qed
\end{thm}
 In fact we have
$$\mathcal{P}_{\mathbb{Z}}^+(1)=\{\boldsymbol{\pi}_{\kappa,J}: \kappa \text{ height function}, \text{ $J\subseteq [1,n]$ interval}\}.$$ The only non-trivial direction is derived as follows. Let $\boldsymbol{\pi}=\boldsymbol{\varpi}_{i_1,a_1}\cdots \boldsymbol{\varpi}_{i_k,a_k}\in \mathcal{P}_{\mathbb{Z}}^+(1)$ such that $a_2=a_1q^{(i_2-i_1+2)}$, then we choose $J=[i_1,i_k]$ and $\kappa$ to be the height function given by $a_{1}=q^{\kappa(i_1)}$ and $i_1,i_3,\dots$ are the sinks and $i_2,i_4,\dots$ are the sources. If  $a_2=a_1q^{-(i_2-i_1+2)}$ we simply change the role of sinks and sources. Then we have $\boldsymbol{\pi}=\boldsymbol{\pi}_{\kappa,J}$.

\begin{rem} One of the main results of \cite{HL10} shows that the category $\mathcal{C}_{q,\kappa}$ is a monoidal categorification of a cluster algebra $\mathcal{A}$ of type $A_n$ when $\kappa$ induces the sink-source orientation on $Q_{\kappa}$ or $\kappa(i)=i$ for all $1\leq i\leq n$. This was later extended by representation theoretic methods to any hight function in \cite{BC19a}. The isomorphism identifies the cluster variables in $\mathcal{A}$ with the prime irreducible objects in $\mathcal{C}_{q,\kappa}$ and cluster monomials are mapped to an equivalence class of an irreducible object in $\mathcal{C}_{q,\kappa}$.
\end{rem}
\section{The main results}\label{section3}
Our main result describes the graded decomposition of certain objects including the graded limits of the prime irreducible objects in the HL category. The main result is formulated in Theorem~\ref{mainthm} for $L(\boldsymbol{\pi})$, $\boldsymbol{\pi}\in \mathcal{P}_{\mathbb{Z}}^+(1)$.
\subsection{} Given a $\mathbb{Z}$-graded space $V=\bigoplus_{r\in\mathbb{Z}}V[r]$, we denote by $\tau_p^{*}V$ the graded vector space whose $r$-th graded piece is $V[r+p]$.
We can write the graded character of $L(\boldsymbol{\pi})$ for $\boldsymbol{\pi}\in\mathcal{P}_{\mathbb{Z}}^+(1)$ as follows
$$\mathrm{ch}(L(\boldsymbol{\pi}))=\sum_{\gamma\in Q^+}[L(\boldsymbol{\pi}): V(\mathrm{wt}(\boldsymbol{\pi})-\gamma)]_q  \ \mathrm{ch}_{\lie h}( V(\mathrm{wt}(\boldsymbol{\pi})-\gamma))$$
where
$$\big[L(\boldsymbol{\pi}): V(\mathrm{wt}(\boldsymbol{\pi})-\gamma)\big]_q= \sum_{p=0}^{\infty}\ [L(\boldsymbol{\pi}): \tau_p^{*}V(\mathrm{wt}(\boldsymbol{\pi})-\gamma)] \cdot q^p.$$
In order to state the results, we need some more notation. Fix a dominant integral weight $\lambda\in P^+$ and a multipartition $$\boldsymbol{\mu}=(\mu_1,\dots,\mu_n),\ \ |\mu_i|=r_i,\ \ \mu_i=(\mu_i^{1}\geq \mu_i^{2}\geq\dots\geq \mu_i^{r_i}\geq 0).$$ 
We denote the number of boxes in the first $s$ columns of $\mu_i$ by $\mu_i(s)$, the number of rows of length $r$ in $\mu_i$ by $m_{i,r}$ and $d(\mu_i)$ is defined as the total number of rows of $\mu_i$. Set
$$P^{\boldsymbol{\mu},\lambda}_{s,i}=\lambda_i-2\mu_i(s)+\mu_{i-1}(s)+\mu_{i+1}(s),\ \ \ \ \ \ \ \ \  1\leq i\leq n,\ 1\leq s\leq r_i $$
and
$$K^{\lambda}_{\boldsymbol{\mu}}=\sum_{i=1}^n\Big(\sum_{j=1}^{r_i}\big(2j\mu_i^{j}-\mu_{i+1}(\mu_i^{j})\big)-\lambda_i \cdot d(\mu_i)\Big).$$
\begin{example}\label{kmuausg} Let $n=8$ and $\lambda=\varpi_2+\varpi_3+\varpi_4+\varpi_5+\varpi_7$. 
For the multipartitions 

$$\boldsymbol{\mu}_1={\fontsize{6}{6}\selectfont \Yvcentermath1	
			\Big(\yng(1)\,,\yng(2,1)\,,\yng(2,2)\,,\yng(2,2)\,,\yng(2,1)\,,\yng(2)\,,\yng(1)\,,\emptyset\Big)},\quad \boldsymbol{\mu}_2={\fontsize{6}{6}\selectfont \Yvcentermath1	\Big(\yng(1)\,,\yng(2,1)\,,\yng(2,1,1)\,,\yng(2,1,1)\,,\yng(1,1,1)\,,\yng(1,1)\,,\yng(1)\,,\emptyset\Big)}$$
			
a direct calculation gives $K^{\lambda}_{\boldsymbol{\mu}_1}=13$ and $K^{\lambda}_{\boldsymbol{\mu}_2}=10$.
\end{example}

We will prove the following result in the rest of the paper.
\begin{thm}\label{mainthm} Let $\gamma=\sum_{i=1}^n r_i \alpha_i\in Q^+$, $p\in\mathbb{Z}_+$ and $\boldsymbol{\pi}=\boldsymbol{\varpi}_{i_1,a_1}\cdots \boldsymbol{\varpi}_{i_k,a_k}\in \mathcal{P}_{\mathbb{Z}}^+(1)$. Then 
$$[L(\boldsymbol{\pi}): \tau_p^{*}V(\mathrm{wt}(\boldsymbol{\pi})-\gamma)]=\sum_{\boldsymbol{\mu}} L^p_{\boldsymbol{\mu}}$$ where the sum runs over all mutipartitions $\boldsymbol{\mu}=(\mu_1,\dots,\mu_n)$ with $|\mu_i|=r_i$ and $L^p_{\boldsymbol{\mu}}$ is the number of lattice points in the polytope consisting of all points $(C_{d,r,i})$ $(1\leq i\leq n, 1\leq r\leq r_i\
, 1\leq d\leq m_{i,r})$ satisfying the following inequalities:
$$C_{d,r,i}\geq 0,\ \forall d,r,i,\ \ \sum_{d=1}^{m_{i,r}}C_{d,r,i}\leq P^{\boldsymbol{\mu},\mathrm{wt}(\boldsymbol{\pi})}_{r,i} \ \forall r,i$$
$$\sum_{i=i_j}^{i_{j+1}}C_{m_{i,1},1,i}\geq 1,\ \forall j\in[1,k),\ \ \sum_{i,r,d}d\cdot C_{d,r,i}=-p+|\gamma|-K^{\mathrm{wt}(\boldsymbol{\pi})}_{\boldsymbol{\mu}}.$$
\end{thm}
\subsection{}
In some cases the graded limit remains irreducible. For example, if we have $\boldsymbol{\pi}=\boldsymbol{\varpi}_{i_1,a}\boldsymbol{\varpi}_{i_2,aq^{\pm(i_2-i_1+2)}}$ for some $a\in\mathbb{C}(q)$ and $1\leq i_1< i_2\leq n$, then
\begin{equation}\label{irr1a}V_q(\boldsymbol{\pi})\cong V_q(\varpi_{i_1}+\varpi_{i_2})\end{equation}
and hence we need to have 
$$|L_{\boldsymbol{\mu}}|=0,\ \text{ where } L_{\boldsymbol{\mu}}=\displaystyle\sum_{p=0}^{\infty} L^p_{\boldsymbol{\mu}},$$
unless $\gamma=0$ and $\boldsymbol{\mu}$ is the tuple of empty partitions. To see this from Theorem~\ref{mainthm}, we fix $\gamma\neq 0$ and a non-trivial multipartition $\boldsymbol{\mu}$ and suppose that $u\in\{1,\dots,n\}$ is maximal such that $r_u\neq 0$. Let $\mathrm{wt}(\boldsymbol{\pi})=\lambda=\varpi_{i_1}+\varpi_{i_2}$ and note that $|L_{\boldsymbol{\mu}}|\neq 0$ only if we have $P^{\boldsymbol{\mu},\lambda}_{r_i,i} \geq 0$. This gives
$$2r_i\leq r_{i-1}+r_{i+1}+\lambda_i$$
and hence
$$(i+1)r_i\leq i r_{i+1}+\lambda_1+2\lambda_2+\cdots+i\lambda_i.$$
Now adding $(i+1)\lambda_{i+1}$ on both sides of the inequality and dividing by $(i+1)$ gives
\begin{equation}\label{tttz1}r_i+\lambda_{i+1}\leq  r_{i+1}+1,\ \ 1\leq i\leq n.\end{equation}
Moreover, another neccessary condition of $L_{\boldsymbol{\mu}}$ beeing non-empty is coming from $P^{\boldsymbol{\mu},\lambda}_{1,i}\geq 0$:
\begin{equation}\label{tttz2}\mu_{i-1}^{1}\geq 2\mu_{i}^{1}-\mu_{i+1}^{1}-\lambda_i,\ \ 1\leq i\leq n.\end{equation}
So we have the following consequences from \eqref{tttz1} and \eqref{tttz2}:
\begin{enumerate}
\item [(a)] $\mu_u$ is a single box \vspace{0,05cm} 
\item [(b)] If $r_{i+1}=r_{i+2}+1,$ then $r_{i}=r_{i+1}-\lambda_{i+1}+1$. \vspace{0,05cm} 
\item  [(c)] If $r_{i+1}=r_{i+2}$, then $r_i=r_{i+1}-\lambda_{i+1}$  or $r_i=r_{i+1}-\lambda_{i+1}+1$.\vspace{0,05cm} 
\item [(d)] If $r_{i+1}=r_{i+2}-1$, then $r_i=r_{i+1}-\lambda_{i+1}-1$ or $r_i=r_{i+1}-\lambda_{i+1}+1$ or $r_i=r_{i+1}-\lambda_{i+1}$ .
\end{enumerate}

\vspace{0,1cm}

\textit{Case 1:} In this case we suppose that $u<i_1$. Since $r_u=1=r_{u+1}+1$ (see part (a)) we obtain from part (b) that $r_{u-1}=2$ and from \eqref{tttz2} that $\mu_{u-1}$ is a column tableaux. Continuing in the same fashion we get that $r_i=u-i+1$ and $\mu_i$ is a column tableaux for all $i\in\{1,\dots,u\}$. But this contradicts $2r_1=2u\leq r_2=u-1$ which we have from \eqref{tttz2}.

\vspace{0,1cm}

\textit{Case 2:} In this case we suppose that $u\geq i_2$. Similarly as in the above case we can show that
$$r_i=u-i+1-\lambda_{i+1},\ \text{and}\ \mu_i \text{ is a column tableaux for $i_2-1\leq i\leq u$}.$$
In particular, $r_{i_2-1}=r_{i_2}$, which means that we have two choices for $r_{i_2-2}$ by part (c). Either $r_{i_2-2}=r_{i_2-1}$ or $r_{i_2-2}=r_{i_2-1}+1$. But part (b) forces in the latter case that we have to keep increasing until we reach index $i_1$:
$$r_{i_1-1}=r_{i_1}>\cdots >r_{i_2-2}>r_{i_2-1}.$$
Again by part (c) we have $r_{i_1-2}=r_{i_1-1}$ or $r_{i_1-2}=r_{i_1-1}+1$ and hence we have a weakly increasing sequence
$$r_1\geq r_2\geq \cdots \geq r_{i_1-1}$$
which contradicts once more $2r_1\leq r_2$. In conclusion we must have $r_{i_2-2}=r_{i_2-1}=r_{i_2}$. Continuing with the same argument, we get
$$r_{i_1-1}+1=r_{i_1}=\cdots=r_{i_2}=r_{i_2+1}+1$$
and each partition $\mu_{i_1},\dots,\mu_{i_2}$ is a column tableaux by \eqref{tttz2}. This gives $P^{\boldsymbol{\mu},\lambda}_{1,i}=0$ for $i\in\{i_1,\dots,i_2\}$ and hence $C_{r_{i_1},1,i_1}+\cdots+C_{r_{i_2},1,i_2}=0$ which is a contradiction.

\vspace{0,1cm}

\textit{Case 3:} The last case $i_1\leq u<i_2$ works similar and we omit the details.

This proves that $L_{\boldsymbol{\mu}}$ can only be non-empty if $\boldsymbol{\mu}$ consists of empty-partitions and we get \eqref{irr1a} from Theorem~\ref{mainthm}.

\begin{example} Let $n=8$ and set 
$$\gamma=\alpha_1+3\alpha_2+4\alpha_3+4\alpha_4+3\alpha_5+2\alpha_6+\alpha_7,\ \ i_1=2,\ i_2=3,\ i_3=4,\ i_4=5,\ i_5=7.$$
Note that only multipartitions $\boldsymbol{\mu}$,
for which $P^{\boldsymbol{\mu},\lambda}_{r,i} \ge 0$ for all $1 \le i \le n$ and $1 \le r \le r_i$, need to be considered. A long but tidious calculation shows that the relevant multipartitions are given by 
$$ \fontsize{6}{6}\selectfont \Yvcentermath1
      \Big(\yng(1)\,,\yng(1,1,1)\,,\yng(1,1,1,1)\,,\yng(1,1,1,1)\,,\yng(1,1,1)\,,\yng(2)\,,\yng(1)\,,\emptyset\Big), 
      \Big(\yng(1)\,,\yng(2,1)\,,\yng(2,1,1)\,,\yng(2,1,1)\,,\yng(2,1)\,,\yng(2)\,,\yng(1)\,,\emptyset\Big), 			
			\Big(\yng(1)\,,\yng(2,1)\,,\yng(2,2)\,,\yng(2,2)\,,\yng(2,1)\,,\yng(2)\,,\yng(1)\,,\emptyset\Big),$$
$$\fontsize{6}{6}\selectfont \Yvcentermath1 \left.
      \Big(\yng(1)\,,\yng(2,1)\,,\yng(2,2)\,,\yng(2,1,1)\,,\yng(1,1,1)\,,\yng(1,1)\,,\yng(1)\,,\emptyset\Big), 
      \Big(\yng(1)\,,\yng(1,1,1)\,,\yng(1,1,1,1)\,,\yng(1,1,1,1)\,,\yng(1,1,1)\,,\yng(1,1)\,,\yng(1)\,,\emptyset\Big), 
      \Big(\yng(1)\,,\yng(2,1)\,,\yng(2,1,1)\,,\yng(2,1,1)\,,\yng(1,1,1)\,,\yng(1,1)\,,\yng(1)\,,\emptyset\Big)
\right.$$
Moreover, the polytopes described in Theorem~\ref{mainthm} are only non-empty for the multipartitions
$$ \boldsymbol{\mu}_1 = {\fontsize{6}{6}\selectfont \Yvcentermath1 \Big(\yng(1)\,,\yng(2,1)\,,\yng(2,2)\,,\yng(2,2)\,,\yng(2,1)\,,\yng(2)\,,\yng(1)\,,\emptyset\Big)}, \quad 
\boldsymbol{\mu}_2={\fontsize{6}{6}\selectfont \Yvcentermath1 \Big(\yng(1)\,,\yng(2,1)\,,\yng(2,1,1)\,,\yng(2,1,1)\,,\yng(1,1,1)\,,\yng(1,1)\,,\yng(1)\,,\emptyset\Big)}.$$
Using Example~\ref{kmuausg} we obtain that 
\begin{align*}
  L^p_{\boldsymbol{\mu}_1} &=| \lbrace C_{1,2,4} \in \mathbb{Z}_{+} \,\colon\, C_{1,2,5}\le 1 , \,   C_{1,2,5} + p = 5\rbrace| \\
  L^p_{\boldsymbol{\mu}_2} &=| \lbrace (C_{1,1,2}, C_{2,1,4}, C_{1,1,7}) \in \mathbb{Z}_{+}^3 \colon\, 
C_{1,1,2} = C_{2,1,4} = C_{1,1,7} = 1,\, p =4 \rbrace|
\end{align*}
and hence
$$[L(\boldsymbol{\pi}):V(\mathrm{wt}(\boldsymbol{\pi})-\gamma)]_q=  2q^4 + q^5.$$
\end{example}
\section{Truncated modules for current algebras}\label{section4}
\subsection{}
We introduce a class of truncated $\mathfrak{g}[t]$-modules indexed by
pairs $(\boldsymbol{\xi},\lambda)$ where $\boldsymbol{\xi}=(\xi_{\alpha})_{\alpha\in R^+}$ is a tuple of non-negative integers indexed by the set of positive roots of $\mathfrak{g}$ and $\lambda$ is a dominant integral weight (abbreviate as usual $\xi_{\alpha_{i,j}}=\xi_{i,j}$ etc.). We define $M_{\boldsymbol{\xi},\lambda}$ to be the cyclic quotient of the local Weyl module $W(\lambda)$ by the submodule generated by
$$\left\{(x^{-}_{\alpha}\otimes t^{\xi_{\alpha}})w_{\lambda} : \alpha \in R^+\right\}$$
where $w_{\lambda}$ is the highest weight generator of $W(\lambda)$. In other words, $M_{\boldsymbol{\xi},\lambda}\cong \mathbf{U}/\mathcal{I}_{\boldsymbol{\xi},\lambda}$, where $\mathcal{I}_{\boldsymbol{\xi},\lambda}\subseteq \mathbf{U}$ is the left ideal generated by the elements
$$\mathfrak{n}^+[t],\ h_i\otimes t^r,\  h_i-\lambda_i\cdot 1\ (r\in\mathbb{N},\ 1\leq i\leq n),$$
$$(x_i^{-}\otimes 1)^{\lambda_i+1}\ (1\leq i\leq n),\ (x^{-}_{\alpha}\otimes t^{\xi_{\alpha}}),\ \alpha\in R^+.$$ 

We call $\boldsymbol{\xi}$ \textit{normalized} if $\xi_{\alpha}=\min\{\xi_{\beta}: \beta\succeq \alpha\}$ for all $\alpha\in R^+$. Given an arbitrary $\boldsymbol{\xi}$ we obviously have $M_{\boldsymbol{\xi},\lambda}\cong M_{\boldsymbol{\xi}',\lambda}$ where $\boldsymbol{\xi}'$ is defined by $\xi_{\alpha}':=\min\{\xi_{\beta}: \beta\succeq \alpha\}$. Hence we can replace $\boldsymbol{\xi}$ by the normalized tuple $\boldsymbol{\xi}'$ and suppose in the rest of the paper without loss of generality that $\boldsymbol{\xi}$ is normalized unless otherwise stated.
\begin{example}\label{beisp1}
\begin{enumerate}
\item If we choose $\xi_{\alpha}=1$ for all $\alpha\in R^+$ we obtain the pull-back of the irreducible representation $M_{\boldsymbol{\xi},\lambda}=\mathrm{ev}_0^{*}V(\lambda)$ where $\mathrm{ev}_0: \mathfrak{g}[t]\rightarrow \mathfrak{g}$ is the evaluation map at $0$.
\item Let $N\in\mathbb{N}$. Choosing $\xi_{\alpha}=N$ for all $\alpha\in R^+$ we obtain the truncated Weyl module whose structure has been partially determined in the series of articles \cite{BK20,FMM19,KL14}.
\item Let $\boldsymbol{\pi}\in\mathcal{P}^+_{\mathbb{Z}}(1)$. Then $L(\boldsymbol{\pi})$ appears as $M_{\boldsymbol{\xi},\mathrm{wt}(\boldsymbol{\pi})}$ where $\xi_{\alpha}=\lceil \frac{\mathrm{wt}(\boldsymbol{\pi})(h_{\alpha})}{2}\rceil$ for all $\alpha\in R^+$.
\end{enumerate}
In fact many families of Demazure modules and generalized Demazure modules appear as $M_{\boldsymbol{\xi},\lambda}$ for a suitable choice of $\boldsymbol{\xi}$ and $\lambda$; see for example \cite{CV13,KV14} and \cite{CDM19}.
\end{example}
\subsection{} The aim of this section is to describe the $\mathbf{U}^{-}$ structure of $M_{\boldsymbol{\xi},\lambda}$. We set $x^{(r)}:=\frac{1}{r!} x^r$ for an element $x\in \lie g[t]$, $r\in\mathbb{Z}_+$ and $$\mathbf{x}_i^{-}(r,s):=\sum_{} (x_i^{-}\otimes 1)^{(b_0)}\cdots (x_i^{-}\otimes t^s)^{(b_s)}$$
where the sum runs over all tuples $(b_0,\dots,b_s)$ of non-negative integers satisfying $b_0+\cdots+b_s=r$ and $b_1+2b_2+\cdots + sb_s=s$. The following result is a slight modification of \cite[Lemma 20]{BKu17}.
\begin{lem}\label{gar}Let $V$ be a $\mathbf{U}$-representation and $v\in V$ a weight vector such that $$(h\otimes t^{r+1})v=(x_i^+\otimes t^r)v=0,\ \text{for all } h\in \lie h,\ r\in \mathbb{Z}_+.$$ Then we have for all $r,s\in \mathbb{Z}_+$ with $r+s\geq 1+\lambda_i$:
\begin{equation}\label{str}\mathbf{U}\cdot \mathbf{x}_i^{-}(r,s).v \subseteq \sum_{u+w\geq 1+\lambda_i}\mathbf{U}^{-}\cdot \mathbf{x}_i^{-}(u,w)v\end{equation}
\begin{proof}
From \cite[Proof of Lemma 20]{BKu17} we obtain up to a non-zero constant
\begin{align*}(h\otimes t^{\ell})\cdot \mathbf{x}_i^{-}(r,s)v=r\cdot \mathbf{x}_i^{-}(r,s+\ell)v-\sum_{j=0}^{\ell-1}(x_i^{-}\otimes t^j)\cdot \mathbf{x}_i^{-}(r-1,s+\ell-j)v.\end{align*}
Hence the above element lies in the $\mathbf{U}^-$--span of elements of the form $\mathbf{x}_i^{-}(u,w)v$ with $u+w\geq 1+\lambda_i$. The fact that an arbitrary product of elements in $(\mathfrak{h}\otimes t\mathbb{C}[t])$ applied to $\mathbf{x}_i^{-}(r,s)v$ lies in the right hand side of \eqref{str} follows from $[\lie h,\lie n^{-}]\subseteq \lie n^-$ and induction on the length. Now we consider the element $(x_i^+\otimes t^{a})\cdot\mathbf{x}_i^{-}(r,s)v$. If $a=1$, then clearly (up to a constant)
$$(x_i^+\otimes t)\cdot \mathbf{x}_i^{-}(r,s)v=\mathbf{x}_i^{-}(r-1,s+1)v.$$
Otherwise, we choose $h\in\lie h$ with $\alpha_i(h)\neq 0$ and get
$$(x_i^+\otimes t^{a})\cdot \mathbf{x}_i^{-}(r,s)v=\big[(h\otimes t),(x_i^+\otimes t^{a-1})\big]\cdot \mathbf{x}_i^{-}(r,s)v$$
and the claim in this case (length one) follows by induction on $a$. Again the general case follows by induction on the length and the Poincaré--Birkhoff--Witt theorem.
\end{proof}
\end{lem}
We have an isomorphism of $\mathbf{U}^{-}$--modules: $$\mathbf{U}^{-}/\mathcal{J}_{\boldsymbol{\xi},\lambda} \xrightarrow{\sim} M_{\boldsymbol{\xi},\lambda},\ \ \mathcal{J}_{\boldsymbol{\xi},\lambda}:=(\mathcal{I}_{\boldsymbol{\xi},\lambda}\cap \mathbf{U}^{-}).$$
\begin{prop}\label{relc}
We have that $\mathcal{J}_{\boldsymbol{\xi},\lambda}$ is the left ideal in $\mathbf{U}^{-}$ generated by the elements 
\begin{align}&\label{1}\mathbf{x}_{i}^{-}(r,s),\ 1\leq i\leq n,\ r\in \mathbb{N},\ s\in\mathbb{Z}_+:\ r+s\geq 1+\lambda_i&\\& \label{2}(x^{-}_{\alpha}\otimes t^{k}),\ \alpha\in R^+,\ k\geq \min\{\xi_{\beta}: \beta\succeq \alpha\}.\end{align}
\begin{proof}
Obviously each element in \eqref{1} and \eqref{2} is contained in $\mathcal{J}_{\boldsymbol{\xi},\lambda}$, where the containment of \eqref{1} follows from the Garland identities (see for example \cite[Lemma 7.1]{G78}) and the containment of \eqref{2} follows from \eqref{seqr}. In order to finish the proof we do the following procedure. We write each element in $\mathcal{I}_{\boldsymbol{\xi},\lambda}$ in PBW order:
$$\mathbf{U}\cong \mathbf{U}\cdot\mathfrak{n}^+[t]\ \oplus \ \mathbf{U}^{-}\cdot \mathbf{U}^0\cdot \mathfrak{h}[t]_+\ \oplus \ \mathbf{U}^{-}\cdot\mathbf{U}(\mathfrak{h}).$$
After that, we take the projection of that element onto $\mathbf{U}^{-}\cdot\mathbf{U}(\mathfrak{h})$ and substitute $h=\lambda(h)$ for all $h\in\lie h$. This gives an element in $\mathcal{J}_{\boldsymbol{\xi},\lambda}$ and it is clear that each element in $\mathcal{J}_{\boldsymbol{\xi},\lambda}$ appears in this way. So the strategy of the proof is to show that by doing this procedure we only get elements in the left ideal in $\mathbf{U}^{-}$ generated by the elements \eqref{1} and \eqref{2}. Note that it is enough to project all elements in
 \begin{equation}\label{3}\sum_{i=1}^n \mathbf{U}\cdot (x_i^-\otimes 1)^{\lambda_i+1}+\sum_{\alpha\in R^+}\mathbf{U}\cdot (x^-_{\alpha}\otimes t^{\xi_{\alpha}}).\end{equation}
Using Lemma~\ref{gar} we see that the projection of each element in \eqref{3} is contained in 
 $$\sum_{i=1}^n\sum_{\substack{r,s:\\ r+s>\lambda_i}} \mathbf{U}^{-}\cdot\mathbf{x}^-_i(r,s)+\sum_{\alpha\in R^+}\sum_{k\geq \min\{\xi_{\beta}: \beta\succeq \alpha\}}\mathbf{U}^{-}\cdot (x^-_{\alpha}\otimes t^{k})\mathbf{U}(\lie h).$$
This finishes the proof of the proposition.
\end{proof}
\end{prop}
\subsection{}Our next goal is to understand the left ideal $\mathcal{L}_{\boldsymbol{\xi},\lambda}:=\mathbf{U}_{\ell}^-\mathcal{J}_{\boldsymbol{\xi},\lambda}$ generated by 
$\mathcal{J}_{\boldsymbol{\xi},\lambda}$ in the loop algebra. The reason is that the same proof as in \cite[Lemma 4.2]{N17a} shows that we have a linear isomorphism 
\begin{equation}\label{ab2}\mathbf{U}_{\ell}^-/\left(\mathfrak{n}^-[t^{-1}]\mathbf{U}_{\ell}^{-}+\mathcal{L}_{\boldsymbol{\xi},\lambda})\right)\cong \mathbf{U}^-/(\mathfrak{n}^-\mathbf{U}^{-}+\mathcal{J}_{\boldsymbol{\xi},\lambda})\end{equation}
preserving the $-Q^+$-grading. This isomorphism is induced from the projection $\mathbf{U}_{\ell}^-\rightarrow \mathbf{U}^-$ with respect to the decomposition
$$\mathbf{U}_{\ell}^-\cong (\mathfrak{n}^-\otimes \mathbb{C}[t^{-1}]_+)\mathbf{U}_{\ell}^-\oplus \mathbf{U}_{}^-$$
and its importance will be discussed later. The explicit generators of $\mathcal{L}_{\boldsymbol{\xi},\lambda}$ are given as follows.
Let $N=\text{max}\{\lambda_i,\xi_{i}: 1\leq i\leq n\}$ and define
\begin{equation*}\widehat{X}_i^{-}(z)=\sum^N_{m=-\infty} (x_{i}^{-}\otimes t^m) z^{-m-1},\ 1\leq i\leq n.\end{equation*}
We denote by $\widehat{\mathbf{x}}_{i}^{-}(r,s)$ the coefficient in front of $z^{-(r+s)}$ in the series $\left(\widehat{X}_i^{-}(z)\right)^r$.

\begin{lem}\label{unminus}The $\mathbf{U}_{\ell}^-$ left ideal $\mathcal{L}_{\boldsymbol{\xi},\lambda}$ is generated by the elements 
\begin{align}&\notag \widehat{\mathbf{x}}_{i}^{-}(r,s),\ 1\leq i\leq n,\ r\in \mathbb{N},s\in\mathbb{Z}_+:\ r+s\geq 1+\lambda_i&\\& \label{gen2}(x^{-}_{\alpha}\otimes t^{k}),\ \alpha\in R^+,\ k\geq \min\{\xi_{\beta}: \beta\succeq \alpha\}.\end{align}

\begin{proof}
The proof is similar to the proof of \cite[Lemma 4.2]{N17a} using Proposition~\ref{relc}.
\end{proof}
\end{lem}
Now we come back to the importance of the graded version of \eqref{ab2}. Recall that $M_{\boldsymbol{\xi},\lambda}[p]$ denotes the $p$-th graded piece of $M_{\boldsymbol{\xi},\lambda}$ for $p\in\mathbb{Z}_+$ (similarly for $\mathbf{U}^{-}$). Setting
$$\mathbf{U}_{\ell}^-[p]:=(\mathfrak{n}^-\otimes \mathbb{C}[t^{-1}]_+)\mathbf{U}_{\ell}^-\oplus \mathbf{U}_{}^-[p]
$$
we obtain 
\begin{align}[M_{\boldsymbol{\xi},\lambda}: \tau^{*}_pV(\lambda-\gamma)]\notag&=\text{dim} \left(M_{\boldsymbol{\xi},\lambda}[p]/\mathfrak{n}^-M_{\boldsymbol{\xi},\lambda}[p]\right)_{\lambda-\gamma}&\\&=\notag\text{dim}\left(\mathbf{U}^-[p]/\left(\mathfrak{n}^-\mathbf{U}^-[p]+\mathcal{J}_{\boldsymbol{\xi},\lambda}\cap \mathbf{U}^-[p]\right)\right)_{-\gamma}&\\&\label{ab1}=\text{dim} \left(\mathbf{U}_{\ell}^-[p]/\left((\mathfrak{n}^-\otimes \mathbb{C}[t^{-1}]_+)\mathbf{U}_{\ell}^-+\mathfrak{n}^-\mathbf{U}^-[p]+\mathcal{L}_{\boldsymbol{\xi},\lambda}\cap\mathbf{U}_{\ell}^-[p]\right)\right)_{-\gamma}\end{align} Hence the graded character of $M_{\boldsymbol{\xi},\lambda}$ is determined by the numbers \eqref{ab1}.
\section{Dual functional realization of loop algebras}\label{section5} In this subsection we review the dual functional realization of $\mathbf{U}_{\ell}^-$ from \cite{FS94}; see also \cite{AK07,AKS06,FKLMM02} for further developments.

\subsection{} Let $\gamma = \sum^n_{i=1}r_i \alpha_i$ be an element of $Q^+$. Consider the generating current
$$x_\alpha^{-}(z)=\sum_{r\in\mathbb{Z}} (x_{\alpha}^{-}\otimes t^r) z^{-r-1},\ \alpha\in R^+.$$
The graded component $\left(\mathbf{U}_{\ell}^-\right)_{-\gamma}$ with respect to the $\mathfrak{h}$-grading is generated by the coefficients of products of generating currents of the form
$$x^{-}_{i_1}(z_1)\cdots x^{-}_{i_{k}}(z_{k}),\ \gamma=\alpha_{i_1}+\alpha_{i_2}+\cdots+\alpha_{i_k}.$$ 
Note that we have two types of operator product expansion relations induced from the commutation relations and the Serre relations (see for example \cite[Section 4]{AK07}). Define the formal delta function by $\delta(z_1-z_2)=\sum_{r\geq 0}z_1^{-r-1}z_2^{r}$ and for a series
$f(z)=\sum_{r\in\mathbb{Z}} f_r z^{-r-1}$ we set
$$f(z)_{-}=\sum_{r\geq 0} f_r z^{-r-1},\ f(z)_{+}=\sum_{r<0} f_r z^{-r-1}.$$
\subsection{} We will see at the end of this subsection that a typical element in the dual space can be viewed as a rational function in the variables $$\bold{x}_{\gamma}=\{ x_{i,r} : 1\leq i \leq n , 1 \leq r \leq r_i\}.$$ 
Due to the restrictions coming from the OPE relations we will get some restrictions on the rational functions. Let us consider products of the form 
$$x_{\alpha_{i_1}}^{-}(x_{i_1,r_1})\cdots x_{\alpha_{i_s}}^{-}(x_{i_s,r_s})=\sum_{k_1,\dots,k_s\in\mathbb{Z}}(x_{\alpha_{i_1}}\otimes t^{k_1})\cdots (x_{\alpha_{i_s}}\otimes t^{k_s})x_{i_1,r_1}^{-k_1-1}\cdots x_{i_s,r_s}^{-k_s-1}$$
so that a Laurant series in the variables $\bold{x}_{\gamma}$ can be viewed naively as an element in the dual space via the above product. The first OPE relation  (see \cite[Section 4]{AK07}) states that the functions will have at most a simple pole whenever $x_{i,r}=x_{i+1,s}$. The Serre relations tell us that the evaluation of the functions at $x_{i,1}=x_{i,2}=x_{i\pm 1,1}$ must vanish. Moreover, since $[x_{\alpha_{i}}^{-}(z),x_{\alpha_{i}}^{-}(w)]=0$ the functions have a symmetry under the exchange of variables $x_{i,r}\leftrightarrow x_{i,s}$. This motivates to 
define $\mathbb{U}_{\gamma}$ to be the subspace of rational functions in the variables $\bold{x}_{\gamma}$ which are of the form $$g(\bold{x}_{\gamma})=\frac{f(\bold{x}_{\gamma})}{\Delta_{\gamma}},\ \ \ \ \Delta_{\gamma}:=\prod_{r,s,i}(x_{i,r}-x_{i+1,s})$$
where $f(\bold{x}_{\gamma})$ is a Laurent polynomial in the variables $\bold{x}_{\gamma}$, symmetric under the action of the parabolic subgroup $S_{r_1}\times \cdots \times S_{r_n}$ ($S_{r_i}$ permutes the variables $x_{i,r},\ 1\leq r\leq r_i$), and vanishing under the specialization $x_{i,1}=x_{i,2}=x_{i\pm 1,1}$. The residue of $g(\bold{x}_{\gamma})\in \mathbb{U}_{\gamma}$ viewed as a function in $x_{i,1}$ is defined as follows. We consider the Laurant series expansion of $g(\bold{x}_{\gamma})$ in a punctured disk $\{0<|x_{i,1}|<\epsilon\}$ by expanding all $(x_{i,1}-x_{i+1,s})^{-1}$, (resp. $(x_{i-1,r}-x_{i,1})^{-1}$) using  
$$-\delta(x_{i+1,s}-x_{i,1})\ \ (\text{resp. } \delta(x_{i-1,r}-x_{i,1})).$$
Then, $\text{Res}_{x_{i,1}}\left(g(\bold{x}_{\gamma})\right)$ is defined to be the coefficient of $\left(x_{i,1}\right)^{-1}$ in the series expansion. We set $$R_{i,p}(g(\bold{x}_{\gamma})):=\text{Res}_{x_{i,1}}\left(\left(x_{i,1}\right)^{p}g(\bold{x}_{\gamma})\right)$$ and view it as a function in $\mathbb{U}_{\gamma-\alpha_i}$ by re-indexing the set $\{ x_{i,2},\dots,x_{i,r_i}\}$ to $\{ x_{i,1},\dots,x_{i,r_i-1}\}$ which is possible by the symmetry of these functions. 
The following theorem can be found in \cite[Theorem 3.3]{AKS06}.
\begin{prop}\label{dualfunctionalrealization}
The space of functions $\mathbb{U}_{\gamma}$ is dual to the graded component $\left(\mathbf{U}^-_{\ell}\right)_{-\gamma}$ with pairing
$$\langle \cdot , \cdot \rangle: \left(\mathbf{U}^-_{\ell}\right)_{-\gamma} \times \mathbb{U}_{\gamma} \rightarrow \mathbb{C}$$ given by the rule \begin{equation}\label{pairing}\langle (x^{-}_{i_1}\otimes t^{k_1})(x^{-}_{i_2}\otimes t^{k_2})\cdots (x^{-}_{i_d}\otimes t^{k_d}),g(\bold{x}_{\gamma})\rangle =R_{i_1,k_1}R_{i_2,k_2}\cdots R_{i_d,k_d}\left(g(\bold{x}_{\gamma})\right).\end{equation}
\hfill\qed
\end{prop}
Using Proposition~\ref{dualfunctionalrealization} we immediately get that the dual space of $$\left(\mathbf{U}_{\ell}^-[p]/\left((\mathfrak{n}^-\otimes \mathbb{C}[t^{-1}]_+)\mathbf{U}_{\ell}^-+\mathfrak{n}^-\mathbf{U}^-[p]+\mathcal{L}_{
\boldsymbol{\xi},\lambda}\cap\mathbf{U}_{\ell}^-[p]\right)\right)_{-\gamma}$$ consists of all functions $g(\mathbf{x}_{\gamma})\in\mathbb{U}_{\gamma}$ satisfying:
\begin{equation}
\label{dual1}
\langle (x_{i}^{-}\otimes t^k)\cdot \left(\mathbf{U}^-_{\ell}\right)_{-\gamma+\alpha_i}, g(\mathbf{x}_{\gamma})\rangle=0,\ \ \forall i\in \{1,\dots,n\},\ k\leq 0
\end{equation}
\begin{equation}
\label{dual0}
\langle \left(\mathbf{U}^-[q]\right)_{-\gamma}, g(\mathbf{x}_{\gamma})\rangle=0,\ \ \forall q\neq p
\end{equation}
\begin{equation}\label{dual2}
    \langle \mathcal{L}_{
\boldsymbol{\xi},\lambda}\cap \left(\mathbf{U}^-_{\ell}\right)_{-\gamma}, g(\mathbf{x}_{\gamma})\rangle=0
\end{equation}
To see this, we only have to show that a function satisfying \eqref{dual1}-\eqref{dual0} and 
$$\langle \mathcal{L}_{
\boldsymbol{\xi},\lambda}\cap \left(\mathbf{U}^-_{\ell}[p]\right)_{-\gamma}, g(\mathbf{x}_{\gamma})\rangle=0$$ also satisfies \eqref{dual2}. This claim is only non-trivial for elements in $\left(\mathbf{U}^-[q]\right)_{-\gamma+r\alpha_i}\widehat{\mathbf{x}}_i^{-}(r,s)$ considered as a subset of $\mathcal{L}_{
\boldsymbol{\xi},\lambda}\cap \left(\mathbf{U}^-_{\ell}\right)_{-\gamma}$. Obviously 
$$\left(\mathbf{U}^-[q]\right)_{-\gamma+r\alpha_i}\widehat{\mathbf{x}}_i^{-}(r,s)\subseteq \left((\mathfrak{n}^-\otimes \mathbb{C}[t^{-1}]_+)\mathbf{U}_{\ell}^-\oplus \mathbf{U}^-[q+s]\right)_{-\gamma}$$
So if $q+s=p$, then the latter element is contained in $\mathcal{L}_{
\boldsymbol{\xi},\lambda}\cap \left(\mathbf{U}^-_{\ell}[p]\right)_{-\gamma}$ and the claim follows by our assumption. If $q+s\neq p$, then the claim follows from \eqref{dual1}-\eqref{dual0}. 

The characterization of functions satisfying \eqref{dual1} is a slight modification of the result \cite[Lemma 4.6]{N17a} and is stated in the next lemma.
\begin{lem}\label{dualprop1}A function $g(\mathbf{x}_{\gamma})=\frac{f(\mathbf{x}_{\gamma})}{\Delta_{\gamma}}\in\mathbb{U}_{\gamma}$ satisfies $\eqref{dual1}$ if and only if
$$\text{deg}_{x_{i,1}}f(\mathbf{x}_{\gamma}) \leq r_{i-1}+ r_{i+1}-2 \text{ for all } i \in \{1,\dots,n\}.$$
\hfill\qed
\end{lem}
\subsection{} In this subsection we investigate the meaning of \eqref{dual0}. Let $\gamma = \sum^n_{i=1}r_i \alpha_i$ be an element of $Q^+$ and set $e_{\gamma}=\sum_{i=1}^{n-1} r_ir_{i+1}$.
\begin{lem}\label{meandual0}
A function $g(\mathbf{x}_{\gamma})=\frac{f(\mathbf{x}_{\gamma})}{\Delta_{\gamma}}\in\mathbb{U}_{\gamma}$ satisfying $\eqref{dual1}$ has in addition property \eqref{dual0}  if and only if $f(\mathbf{x}_{\gamma})$ is homogeneous of degree $-p-|\gamma|+e_{\gamma}$.
\begin{proof}
The property \eqref{dual0} simply means that $g(\mathbf{x}_{\gamma})_{-}$ (terms with strictly negative power with respect to each variable) is homogeneous of degree $-p-|\gamma|$. Since $g(\mathbf{x}_{\gamma})$ satisfies \eqref{dual1}, hence 
$\text{deg}_{x_{i,r}}g(\mathbf{x}_{\gamma}) \leq -2$, we have that $g(\mathbf{x}_{\gamma})$ is homogeneous of degree $-p-|\gamma|$ and the rest follows by expanding each $(x_{j,r}-x_{j+1,s})$ using the $\delta$-function.
\end{proof}
\end{lem}
\subsection{}
It remains to observe the relations coming from \eqref{dual2}. Recall that $(x^{-}_{i}\otimes t^{k})\in \mathcal{L}_{
\boldsymbol{\xi},\lambda}$ for all $k\geq N$. For $r\leq r_i$ we have
$$\left\langle \widehat{X}_i^{-}(z)^r, g(\mathbf{x}_{r \alpha_i})\right\rangle=g(\mathbf{x}_{r \alpha_i})\Big|_{x_{i,1}=\cdots=x_{i,r}=z}$$
and hence the condition
$$\left\langle \left(\mathbf{U}_{\ell}^-\right)_{\gamma-r\alpha_i} \widehat{\mathbf{x}}_i^{-}(r,s), g(\mathbf{x}_{\gamma})\right\rangle=0,\ \ r+s\geq 1+\lambda_i$$
means that the order of the pole of $g(\mathbf{x}_{\gamma})|_{x_{i,1}=\cdots=x_{i,r}=z}$ at $z=0$ can be at most $\lambda_i$. Now we define the main object of this section. 
\begin{defn}\label{mainobj}
Let $
\boldsymbol{\xi}=(\xi_{\alpha})_{\alpha\in R^+}$ be a normalized tuple of positive integers, $p\in\mathbb{Z}_+$ and $\gamma=\sum_{i=1}^nr_i\alpha_i$ be an element of the root lattice. We denote by $\mathcal{V}_{
\boldsymbol{\xi},\gamma,p}$ the subspace of $\mathbb{U}_{\gamma}$ consisting of all homogeneous Laurant polynomials $f(\mathbf{x}_{\gamma})$ of degree $-p-|\gamma|+e_{\gamma}$ satisfying the following properties:
\begin{enumerate}
\item For all $1\leq i\leq n$ we have $\text{deg}_{x_{i,1}}f(\mathbf{x}_{\gamma}) \leq r_{i-1}+r_{i+1}-2$\vspace{0,15cm}

\item For all $i\in\{1,\dots,n\}$ and $1\leq r\leq r_i$ we have 
$$z^{\lambda_i}\cdot f(\bold{x}_{\gamma})\Big|_{x_{i,1}=x_{i,2}=\cdots=x_{i,r}=z}\in\mathbb{C}[x^{\pm}_{i,r}] [z]$$ 
\item For all $1\leq i\leq j\leq n$ we have that
$$z^{\xi_{i,j}}\cdot f(\bold{x}_{\gamma})\Big|_{x_{i,1}=x_{i+1,1}=\cdots=x_{j,1}=z}\in \mathbb{C}[x^{\pm}_{i,r}] [z]$$ 
\end{enumerate}
\end{defn}
Now we state the main theorem of this section.
\begin{thm}\label{wicht} Let $
\boldsymbol{\xi},\gamma,p$ be as in Definition~\ref{mainobj}. We have an isomorphism of vector spaces 
$$\mathcal{V}_{
\boldsymbol{\xi},\gamma,p}\cong \left(\mathbf{U}^-[p]/(\mathfrak{n}^-\mathbf{U}^-[p]+\mathcal{J}_{
\boldsymbol{\xi},\lambda}\cap \mathbf{U}^-[p] )\right)_{-\gamma}^*$$
and hence the graded character of $M_{
\boldsymbol{\xi},\lambda}$ is given by
$$\mathrm{ch} (M_{
\boldsymbol{\xi},\lambda})=\sum_{\gamma\in Q^+} \left(\sum_{p\in\mathbb{Z}_+} \dim(\mathcal{V}_{
\boldsymbol{\xi},\gamma,p})  \ q^p\right)\mathrm{ch}_{\lie h} (V(\lambda-\gamma))$$
\begin{proof} It is clear that the space of all functions in $\mathbb{U}_{\gamma}$ satisfying the properties \eqref{dual1}-\eqref{dual2} describes the aforementioned dual space as a vector space. Using Lemma~\ref{dualprop1} and Lemma~\ref{meandual0} and the discussion preceeding Definition~\ref{mainobj}, it has only to be checked that Definition~\ref{mainobj}(3) is equivalent 
to the fact that
$$\left\langle \left(\mathbf{U}_{\ell}^{-}\right)_{\gamma-\alpha}(x^{-}_{\alpha}\otimes t^k),g(\mathbf{x}_\gamma)\right\rangle =0.$$ But if we write 
$$(x^{-}_{i,j}\otimes t^k)=[(x^{-}_{j}\otimes t^{k_j}),[(x^{-}_{j-1}\otimes t^{k_{j-1}}),\dots,[(x^{-}_{i+1}\otimes t^{k_{i+1}}),(x^{-}_{i}\otimes t^{k_i})]\dots]]$$
for some $k_i,\dots,k_{j}\in \bz$ with $k_1+\dots+k_{\ell}=k$
this follows from 
\begin{align}\label{schw}
   [R_{j,k_j},[R_{i_{j-1},k_{j-1}},&\dots,[R_{i+1,k_{i+1}},R_{i,k_i}]\dots]] g(\mathbf{x}_{\gamma}) &\\&\notag
   =\text{Res}_{x_{j,1}} \Big\{(x_{i,1}-x_{i+1,1})\cdots (x_{j-1,1}-x_{j,1}) g(\mathbf{x}_{\alpha})\Big|_{x_{i,1}=\cdots=x_{j,1}=z}\cdot (x_{j,1})^k\Big\},
\end{align}
which has been proved in \cite[Section 4.2]{N17a}. So we obtain that the order of the pole of 
$$\left((x_{i,1}-x_{i+1,1})\cdots (x_{j-1,1}-x_{j,1})g(\mathbf{x}_{\gamma})\right)|_{x_{i,1}=x_{i+1,1}=\cdots=x_{j,1}=z}$$ at $z=0$ is at most $\xi_{i,j}$.
\end{proof}

\end{thm}

\begin{rem}
A similar statement as above can be proven for any simple finite-dimensional Lie algebra with slightly more complicated property \eqref{schw}. The simplification of property \eqref{schw} to a specialization property as in Definition~\ref{mainobj}(3) also works at least in simply-laced type. Nevertheless we decided to restrict ourself to $\mathfrak{sl}_{n+1}$ (in particular the prime representations in $\mathcal{C}_{q,\kappa}$) since the computation of $\dim V_{
\boldsymbol{\xi},\gamma,p}$ seems to be quite difficult in general. 
\end{rem}

\begin{example}\label{bei2}
Let $n=2$, $\xi_{\alpha}=2$ for each $\alpha\in R^+$, $\gamma=2\alpha_1+\alpha_2$ and $\lambda=7\varpi_1+5\varpi_2$. We consider Laurant polynomials $f(\bold{x}_{\gamma})$ in the variables $\{x_{1,1}, x_{1,2}, x_{2,1}\}$ which are symmetric under the exchange of variables $x_{1,1} \leftrightarrow x_{1,2}$ and satisfy
$$x_{1,1}^{2} f(\bold{x}_{\gamma})\in \mathbb{C}[x_{1,1}, x^{\pm}_{1,2},x^{\pm}_{2,1}],\ \ x_{2,1}^{2} f(\bold{x}_{\gamma})\in \mathbb{C}[x^{\pm}_{1,1}, x^{\pm}_{1,2},x_{2,1}],\ \ \text{deg}_{x_{2,1}}f(\bold{x}_{\gamma})\leq 0$$
$$\text{deg}_{x_{1,1}}f(\bold{x}_{\gamma})\leq -1,\ \ z^{2} f(\bold{x}_{\gamma})\big|_{x_{1,1}=x_{2,1}=z}\in \mathbb{C}[x^{\pm}_{1,2}, z],\ \ f(\bold{x}_{\gamma})\big|_{x_{1,1}=x_{1,2}=x_{2,1}}=0.$$
This is a vector space of dimension 2 with basis
$$f_1(\bold{x}_{\gamma})=\left(x_{1,1}\right)^{-2}\left(x_{1,2}\right)^{-2}\left(x_{2,1}\right)^{-2}\left(\left(x_{2,1}\right)^2+x_{1,1}x_{1,2}-x_{1,2}x_{2,1}-x_{1,1}x_{2,1}\right)$$ $$f_2(\bold{x}_{\gamma})=\left(x_{1,1}\right)^{-2}\left(x_{1,2}\right)^{-2}\left(x_{2,1}\right)^{-2}\left(x_{1,2}\left(x_{2,1}\right)^2+x_{1,1}\left(x_{2,1}\right)^2-2x_{1,1}x_{1,2}x_{2,1}\right).$$
So Theorem~\ref{wicht} implies $[M_{\boldsymbol{\xi},\lambda}: V(\lambda-\gamma)]_q=q^2+q^3$.
\end{example}
\section{Classical decompositions of the prime irreducible objects}\label{section6}
In this section we will focus on a class of representations $M_{\boldsymbol{\xi},\lambda}$ including the modules $L(\boldsymbol{\pi})$ for $\boldsymbol{\pi}\in\mathcal{P}^+_{\mathbb{Z}}(1)$. If we write $\mathrm{wt}(\boldsymbol{\pi})=\omega_{i_1}+\cdots+\omega_{i_k}$, then obviously it is enough to require the last condition in Definition~\ref{mainobj} only for pairs of the form 
 $(i_j,i_{j+1})$ with $1\leq j\leq k-1$. To see this, we simply write 
\begin{align*}\hspace{0,1cm} (x&_{i,j}^{-}\otimes t^{\lceil \frac{m-r}{2}\rceil})=&\\& \big[\cdots\big[\big[\big[(x_{i,i_r-1}^{-}\otimes 1),(x_{i_r,i_{r+1}}^{-}\otimes t)\big],(x_{i_{r+1}+1,i_{r+2}-1}^{-}\otimes 1)\big],(x_{i_{r+2},i_{r+3}}^{-}\otimes t)\big],\dots,(x_{i_{m-1}+1,j}^{-}\otimes 1)\big]\end{align*}
 where $m$ and $r$ are determined by $i_{r-1}<i\leq i_r$ and $i_{m-1}\leq j< i_m$. So we need \eqref{gen2} in Lemma~\ref{unminus} only for the aforementioned pairs.
\subsection{} We first record the following lemma.
\begin{lem}\label{hilf2l} Let $f(\mathbf{x}_{\gamma})$ be a polynomial in the variables $\mathbf{x}_{\gamma}$ such that Definition~\ref{mainobj}(3) holds for all $(i,j)\in\{(i_j,i_{j+1}):1\leq j\leq k-1\}$ with $\xi_{i_j,i_{j+1}}=1$. Then $f(\mathbf{x}_{\gamma})$ can be expressed in the form 
$$f(\mathbf{x}_{\gamma})=\sum_{\mathbf{c}} x_{i_1,1}^{c_{i_1}}x_{i_1+1,1}^{c_{i_1+1}}\cdots x_{i_{k},1}^{c_{i_{k}}}A_{\mathbf{c}}(\mathbf{x}_{\gamma})$$
where the sum runs over all tuples $\mathbf{c}=(c_{i_1},c_{i_1+1},\dots,c_{i_k})$ of non-negative integers satisfying 
$$c_{i_j}+c_{i_j+1}+\cdots+c_{i_{j+1}}= 1,\ \  1\leq j\leq k-1$$
and $A_{\mathbf{c}}(\mathbf{x}_{\gamma})$ is a polynomial in the variables  
$$\{x_{j,r}: 1\leq j\leq n,\  r\neq 1\}\cup\{x_{i,1}: i\notin[i_1,i_k]\}\cup \bigcup_{j=1}^{k-1} \{x_{i,1}: i\in [i_j,i_{j+1}],\ c_{i_j}+\cdots +c_i= 1\}.$$
\proof We prove the claim by induction on $k$. If $k=2$, we can write $f(\mathbf{x}_{\gamma})$ in the form 
$$f(\mathbf{x}_{\gamma})=x_{i_1,1}A_{1}(\mathbf{x}_{\gamma})+\cdots+x_{i_2,1}A_{i_2-i_1+1}(\mathbf{x}_{\gamma})+h(\mathbf{x}_{\gamma})$$
where $h(\mathbf{x}_{\gamma})$ is a polynomial independent of the variables $x_{i_1,1},\dots,x_{i_2,1}$
(otherwise we factor out summands which depend on these variables) and $A_j(\mathbf{x}_{\gamma})$ with $1< j\leq i_2-i_1+1$ does not depend on the variables $x_{i_1,1},x_{i_1+1,1},\dots, x_{i_1+j-2,1}$. But Definition~\ref{mainobj}(3) forces that  $h(\mathbf{x}_{\gamma})$ vanishes under the specialization $x_{i_1,1}=\cdots=x_{i_2,1}$ and hence $h(\mathbf{x}_{\gamma})=0$. This implies the claim for $k=2$. Now by induction hypothesis we can suppose that
\begin{equation}\label{schritt1}f(\mathbf{x}_{\gamma})=\sum_{\mathbf{c}} x_{i_1,1}^{c_{i_1}}\cdots x_{i_{k-1},1}^{c_{i_{k-1}}}A_{\mathbf{c}}(\mathbf{x}_{\gamma})\end{equation}
where the sum runs over all tuples $\mathbf{c}=(c_{i_1},c_{i_1+1},\dots,c_{i_{k-1}})$ satisfying 
$$c_{i_j}+c_{i_j+1}+\cdots+c_{i_{j+1}}= 1,\ \  1\leq j\leq k-2$$
and $A_{\mathbf{c}}(\mathbf{x}_{\gamma})$ is a polynomial in the variables  
$$\{x_{j,r}: 1\leq j\leq n,\  r\neq 1\}\cup \{x_{i,r}: i\notin[i_1,i_{k-1}]\}\cup \bigcup_{j=1}^{k-2} \{x_{i,r}: i\in [i_j,i_{j+1}],\ c_{i_j}+\cdots +c_i= 1\}.$$
We claim that for each summand in \eqref{schritt1} we have
\begin{equation}\label{zuze}A_{\mathbf{c}}(\mathbf{x}_{\gamma})\Big|_{x_{i_{k-1},1}=\cdots=x_{i_k,1}=z}\div z^{1-c_{i_{k-1}}}\end{equation}
and prove this by downward induction using the lexicographic order on $\mathbb{Z}^{i_{k-1}-i_1+1}_+$. If $c_{i_{k-1}}=1$ there is nothing to show; so let $c_{i_{k-1}}=0$. By induction we assume that the claim is true for all $\mathbf{d} \succ \mathbf{c}$. We split the sum in \eqref{schritt1} into three terms
$$f(\mathbf{x}_{\gamma})=x_{i_1,1}^{c_{i_1}}\cdots x_{i_{k-1},1}^{c_{i_{k-1}}}A_{\mathbf{c}}(\mathbf{x}_{\gamma})+\sum_{\mathbf{d}\prec \mathbf{c}} x_{i_1,1}^{d_{i_1}}\cdots x_{i_{k-1},1}^{d_{i_{k-1}}}A_{\mathbf{d}}(\mathbf{x}_{\gamma})+\sum_{\mathbf{d}\succ \mathbf{c}} x_{i_1,1}^{d_{i_1}}\cdots x_{i_{k-1},1}^{d_{i_{k-1}}}A_{\mathbf{d}}(\mathbf{x}_{\gamma}).$$
By induction we have that 
$$\sum_{\mathbf{d}\succ \mathbf{c}} x_{i_1,1}^{d_{i_1}}\cdots x_{i_{k-1},1}^{d_{i_{k-1}}}A_{\mathbf{d}}(\mathbf{x}_{\gamma})\Big|_{x_{i_{k-1},1}=\cdots=x_{i_k,1}=z}$$
is divisible by $z$. Hence if $\widetilde{A}_{\mathbf{d}}(\mathbf{x}_{\gamma})$ denotes the constant term of $A_{\mathbf{d}}(\mathbf{x}_{\gamma})|_{x_{i_{k-1},1}=\cdots=x_{i_k,1}=z}$ with respect to the variable $z$ we get with Definition~\ref{mainobj}(3)
$$x_{i_1,1}^{c_{i_1}}\cdots x_{i_{k-1},1}^{c_{i_{k-1}}}\widetilde{A}_{\mathbf{c}}(\mathbf{x}_{\gamma})+\sum_{\substack{\mathbf{d}\prec \mathbf{c}\\ d_{i_{k-1}}=0}} x_{i_1,1}^{d_{i_1}}\cdots x_{i_{k-1},1}^{d_{i_{k-1}}}\widetilde{A}_{\mathbf{d}}(\mathbf{x}_{\gamma})=0.$$
This implies that the term $$\sum_{\substack{\mathbf{d}\prec \mathbf{c}\\ d_{i_{k-1}}=0}} x_{i_1,1}^{d_{i_1}}\cdots x_{i_{k-1},1}^{d_{i_{k-1}}}\widetilde{A}_{\mathbf{d}}(\mathbf{x}_{\gamma})$$ is contained in the ideal generated by the monomial $x_{i_1,1}^{c_{i_1}}\cdots x_{i_{k-1},1}^{c_{i_{k-1}}}$. In order to show $\widetilde{A}_{\mathbf{c}}(\mathbf{x}_{\gamma})=0$ it will be enough to fix an element $\mathbf{d}\prec \mathbf{c}$ with $d_{i_{k-1}}=0$ and prove that none of the summands of $x_{i_1,1}^{d_{i_1}}\cdots x_{i_{k-1},1}^{d_{i_{k-1}}}\widetilde{A}_{\mathbf{d}}(\mathbf{x}_{\gamma})$ is contained in the aforementioned ideal.
By the definition of the lexicographic order there exists $\ell(\mathbf{d})\in \{i_{j},\dots,i_{j+1}-1\}$ for some $1\leq j\leq k-2$ such that 
$$c_{r}=d_{r}, \text{ for $r<\ell(\mathbf{d})$ and $1=c_{\ell(\mathbf{d})}>d_{\ell(\mathbf{d})}$=0}.$$ Hence we get
$$d_{i_j}+\dots+d_{\ell(\mathbf{d})}<c_{i_j}+\dots+c_{\ell(\mathbf{d})}\leq c_{i_j}+\dots+c_{i_{j+1}}=1$$ and $x_{\ell(\mathbf{d}),1}$ can not appear in $A_{\mathbf{d}}(\mathbf{x}_{\gamma})$ by assumption. This implies the claim and we get \eqref{zuze}. Now using \eqref{zuze}, we can repeat the base case of the induction for each polynomial $A_{\mathbf{c}}(\mathbf{x}_{\gamma})$ in \eqref{schritt1} and obtain the desired expression. 
\endproof
\end{lem}
It is not difficult to translate the dimension of $\mathcal{V}_{\boldsymbol{\xi},\gamma,p}$ into the language of linear algebra by multiplying the Laurant polynomials in Definition~\ref{mainobj} by the factor $\prod_{i,r}x_{i,r}^{\lambda_i}$ and using Lemma~\ref{hilf2l}. 
 We consider a typical polynomial 
$$\sum_{\mathbf{c}}a_{\mathbf{c}} \prod_{i=1}^n x_{i,1}^{c_{i,1}}\cdots x_{i,r_i}^{c_{i,r_i}},\  \ a_{\mathbf{c}} \in \mathbb{C}.$$
By the symmetry we get $a_{\mathbf{c}}=a_{\sigma(\mathbf{c})}$ for each $\sigma\in S_{r_1}\times\cdots\times S_{r_n}$. Hence the coefficients are determined by complex numbers 
$\left(a_{\boldsymbol{\mu}}\right)_{\boldsymbol{\mu}}$
where $\boldsymbol{\mu}=(\mu_1,\dots,\mu_n)$ is a tuple of partitions and each $\mu_i$ is a partition with at most $r_i$ parts and whose entries are bounded by $r_{i-1}+r_{i+1}-2+\lambda_i$ by Definition~\ref{mainobj}(1). Continuing in this way gives a long list of constraints which is in general hard to calculate. In the next subsection we use a different approach developed in \cite{AK07,AKS06} (used also for example in \cite{Tsy20}).
\subsection{} We recall a filtration on the space of rational
functions $\mathbb{U}_{\gamma}$ from \cite[Section 4.1]{AK07}.
Let $\boldsymbol \mu = (\mu_1,\dots,\mu_n)$ be a
multi-partition such that $|\mu_i|=r_i$ for all $1\leq i\leq n$.
Let $m_{i,r}$ the number of parts of length $r$ in the partition $\mu_i$. Our aim is to define a specialization map $$\varphi_{\boldsymbol \mu}: \mathbb{U}_{\gamma}\rightarrow \mathbb{H}_{\boldsymbol \mu},$$
where $\mathbb{H}_{\boldsymbol \mu}$ is the space of rational functions in the variables
$$\mathbf{y}_{\boldsymbol\mu}=\{y_{i,r,u}: 1\leq i\leq n, 1\leq r\leq r_i, 1\leq u \leq m_{i,r}\}.$$
It means that we have a variable for each row in the multi-partition. The specialization map $\varphi_{\boldsymbol \mu}$ does the following: for each $i\in\{1,\dots,n\}$ we fill the boxes of $\mu_i$ with the variables $x_{i,r}$, $1\leq r\leq r_i$, and specialize all variables in the $u$-th row of length $r$ in $\mu_i$ to $y_{i,r,u}$. By the symmetry this definition is independent of the filling.
\begin{example} Let $n=2$ and $\gamma=2\alpha_1+\alpha_1$. Consider the multi-partition
$$\mu_1=\fontsize{6}{6}\selectfont \Yvcentermath1 \yng(2)\ \ \ \mu_2=\fontsize{6}{6}\selectfont \Yvcentermath1\yng(1)$$
If we apply the specialization map to the functions in Example~\ref{bei2} (ignoring the fact that they are not contained in $\mathbb{U}_{\gamma}$) we get 
$$\varphi_{\boldsymbol{\mu}}(f_1(\mathbf{x}_{\gamma}))=y_{1,2,1}^{-4}y_{2,1,1}^{-2}\left(y_{2,1,1}^2+y_{1,2,1}^2-2y_{2,1,1}y_{1,2,1}\right)$$
$$\varphi_{\boldsymbol{\mu}}(f_2(\mathbf{x}_{\gamma}))=2y_{1,2,1}^{-4}y_{2,1,1}^{-2}\left(y_{1,2,1}y_{2,1,1}^2-y_{1,2,1}^2y_{2,1,1}\right)$$
\end{example} 
Now we define a filtration 
\begin{equation}\label{filt1}\{0\} \subseteq \Gamma_{\boldsymbol\mu_1}\subseteq\cdots
\subseteq \Gamma_{\boldsymbol\mu_t}= \mathbb{U}_{\gamma}\end{equation}
as follows. First we define a lexicographical ordering on the set of multi-partitions. We say $\boldmath \nu\succ\boldsymbol\mu$ if there exists some $j\in \{1,\dots,n\}$ such that
$\nu_d=\mu_d$ for all $d<j$ and
$\nu_j>\mu_j$, where the latter is the usual lexicographical order on partitions which is a total order.
Let
\begin{equation*}
\Gamma_{\boldsymbol\mu} = \underset{\boldsymbol\nu>\boldsymbol\mu}{\bigcap} \ker\varphi_{\boldsymbol \nu},\qquad
\Gamma'_{\boldsymbol\mu} = \underset{\boldsymbol \nu\geq \boldsymbol\mu}{\bigcap} \ker\varphi_{\boldsymbol \nu}\subseteq
\Gamma_{\boldsymbol\mu}. 
\end{equation*}
Then we immediately get
$$\Gamma_{\boldsymbol\mu}\subseteq\Gamma_{\boldsymbol\nu},\ \ \text{if}\ \  \boldsymbol\mu<\boldsymbol\nu $$
and hence we obtain the desired filtration \eqref{filt1} (we have $\boldsymbol\mu_i< \boldsymbol\mu_{i+1}$ for all $i$).  We say 
\begin{equation*}
(r,u)<(s,v) \quad \hbox{ if $r>s$ or if $r=s$ and $u<v$.}
\end{equation*}
Let $$\Omega_{\boldsymbol\mu}=\frac{\displaystyle
\prod_{1\leq i \leq n}\prod_{(r,u)<(s,v)}
(y_{i,r,u}-y_{i,s,v})^{2s} } 
{\displaystyle \prod_{1\leq i\leq n}\ \prod_{r,s,u,v}
(y_{i,r,u} - y_{i+1,s,v})^{\min\{r, s\}}}.$$
We obtain the following important result from \cite[Theorem 3.6]{AKS06}.
\begin{lem}\label{AK}
The induced map $\overline{\varphi}_{\boldsymbol \mu}: \Gamma_{\boldsymbol\mu}/\Gamma^{'}_{\boldsymbol\mu}\rightarrow \mathbb{H}_{\boldsymbol\mu}$ of the specialization map $\varphi_{\boldsymbol\mu}$ is an isomorphism of graded vector spaces onto its image. Moreover, the image  is the space of rational functions of the form
\begin{equation*}\label{}
\Omega_{\boldsymbol\mu}\cdot h(\mathbf{y}_{\boldsymbol\mu})
\end{equation*}
where $h(\mathbf{y}_{\boldsymbol\mu})$ is an arbitrary Laurent polynomial in the variables
$\mathbf{y}_{\boldsymbol\mu}$, symmetric with respect to the exchange of variables
$y_{i,r,u}\leftrightarrow y_{i,r,v}$.
\hfill\qed
\end{lem}
Now \eqref{filt1} induces a filtration
$$\{0\} \subseteq \mathcal{V}_{\boldsymbol \xi,\gamma,p}\cap \Gamma_{\boldsymbol\mu_1}\subseteq\cdots
\subseteq \mathcal{V}_{\boldsymbol\xi,\gamma,p}\cap\Gamma_{\boldsymbol\mu_t}= \mathcal{V}_{\boldsymbol \xi,\gamma,p}$$
and hence 
$$\dim(\mathcal{V}_{\boldsymbol\xi,\gamma,p})=\sum_{\boldsymbol\mu} \dim \left((\mathcal{V}_{\boldsymbol\xi,\gamma,p}\cap \Gamma_{\boldsymbol\mu})/(\mathcal{V}_{\boldsymbol\xi,\gamma,p}\cap \Gamma^{'}_{\boldsymbol\mu})\right)=\sum_{\boldsymbol\mu} \dim (\varphi_{\boldsymbol \mu}(\mathcal{V}_{\boldsymbol\xi,\gamma,p}\cap \Gamma^{}_{\boldsymbol\mu})).$$
\subsection{}

In what follows we will describe the functions in $\varphi_{\boldsymbol \mu}(\mathcal{V}_{\boldsymbol \xi,\gamma,p}\cap \Gamma^{}_{\boldsymbol\mu})$; recall the importance of the space $\mathcal{V}_{\boldsymbol\xi,\gamma,p}$ from Theorem~\ref{wicht}. By Definition~\ref{mainobj} we have a list of restrictions to the function $h(\mathbf{y}_{\boldsymbol\mu})$. The second property implies that the Laurent polynomial $h(\mathbf{y}_{\boldsymbol\mu})$ has a pole at $y_{i,r,u}$ of order at most $\lambda_i$. Hence we can write the Laurant polynomials in the image as 
\begin{equation}\label{w2}\frac{\Omega_{\boldsymbol\mu}}{\displaystyle \prod_{1\leq i\leq n}\prod_{(r,u)}(y_{i,r,u})^{\lambda_i}}\ h_1(\mathbf{y}_{\boldsymbol\mu})\end{equation}
where $h_1(\mathbf{y}_{\boldsymbol\mu})$ is an arbitrary polynomial symmetric under the exchange of variables
$y_{i,r,u}\leftrightarrow y_{i,r,v}$. The first condition in Definition~\ref{mainobj} implies that 
$$\text{deg}_{y_{i,r,u}}\varphi_{\boldsymbol\mu}(g(\mathbf{x}_{\gamma}))\leq -2r.$$
Since the degree of the first term in \eqref{w2} with respect to the variable $y_{i,r,u}$ is given by $-P^{\boldsymbol{\mu},\lambda}_{r,i}-2r$ (recall the definition from Section~\ref{section3})
we obtain that 
$$\text{deg}_{y_{i,r,u}}h_1(\mathbf{y}_{\boldsymbol\mu})\leq P^{\boldsymbol{\mu},\lambda}_{r,i}.$$
Since the specialization map preserves the homogeneous degree, we have also that $\varphi_{\boldsymbol\mu}(g(\mathbf{x}_{\gamma}))$ is homogeneous of degree $-p+|\gamma|$ and hence a straightforward translation into the language of partitions gives that $h_1(\mathbf{y}_{\boldsymbol\mu})$ is homogeneous of degree $-p+|\gamma|-K^{\lambda}_{\boldsymbol \mu}$.
Now we investigate the meaning of last condition. Since the multiplication of $f(\mathbf{x}_{\gamma})$ by $\prod_{i,r}x_{i,r}^{\lambda_i}$ is a polynomial, we can apply Lemma~\ref{hilf2l} and obtain that 
\begin{equation}\label{prop3}\displaystyle
\prod_{1\leq i \leq n}\prod_{(r,u)<(s,v)}
(y_{i,r,u}-y_{i,s,v})^{2s} 
\displaystyle \prod_{1\leq i<n}\ \prod_{r,s,u,v}
(y_{i,r,u} - y_{i+1,s,v})^{rs-\min\{r, s\}}\ h_1(\mathbf{y}_{\boldsymbol\mu})\end{equation}
is contained in the intersection of ideals 
\begin{equation}\label{prop4}\mathcal{Y}:=\bigcap_{\substack{s_{i_1},\dots,s_{i_k}\\ v_{i_1},\dots,v_{i_k}}}\mathcal{Y}_{\substack{s_{i_1},\dots,s_{i_k}\\ v_{i_1},\dots,v_{i_k}}}\end{equation}
where
$$\mathcal{Y}_{\substack{s_{i_1},\dots,s_{i_k}\\ v_{i_1},\dots,v_{i_k}}}:=\Big\langle y_{i_1,s_{i_1},v_{i_1}}^{c_{i_1}}y_{i_1+1,s_{i_1+1},v_{i_1+1}}^{c_{i_1+1}}\cdots y_{i_{k},s_{i_k},v_{i_k}}^{c_{i_{k}}}:\ \ \sum_{p=i_j}^{i_{j+1}} c_{p}
= 1,\ \forall  j\in[1,k)\Big\rangle$$
Clearly, a polynomial $f(\mathbf{y}_{\boldsymbol{\mu}})$ is contained in $\mathcal{Y}_{\substack{s_{i_1},\dots,s_{i_k}\\ v_{i_1},\dots,v_{i_k}}}$ if and only if 
\begin{equation}\label{div34} f(\mathbf{y}_{\boldsymbol{\mu}})|_{y_{i_j,s_{i_j},v_{i_j}}=\cdots =y_{i_{j+1},s_{i_{j+1}},v_{i_{j+1}}}=z}\div z\end{equation} for all $j\in [1,k).$ However, this property is already satisfied by the prefactor of $h_1(\mathbf{y}_{\boldsymbol{\mu}})$ in \eqref{prop3} unless there exists $j\in [1,k)$ with $s_{i_j}=\cdots=s_{i_{j+1}}=1$. Moreover, if there exists such a $j$ with $s_{i_j}=\cdots=s_{i_{j+1}}=1$, then the prefactor is a polynomial which has a non-zero constant term after evaluating the variables $y_{i_j,s_{i_j},v_{i_j}}=\cdots=y_{i_{j+1},s_{i_{j+1}},v_{i_{j+1}}}=z$ (viewed as a polynomial in $z$).  Hence \eqref{prop3} is contained in $\mathcal{Y}_{\substack{s_{i_1},\dots,s_{i_k}\\ v_{i_1},\dots,v_{i_k}}}$ if and only if $h_1(\mathbf{y}_{\boldsymbol{\mu}})$ satisfies \eqref{div34} for all $j\in [1,k)$ with $s_{i_j}=\cdots=s_{i_{j+1}}=1$. This implies the following lemma.
\begin{lem}\label{conta}The containement of \eqref{prop3} in the ideal \eqref{prop4} is equivalent to the fact that $h_1(\mathbf{y}_{\boldsymbol{\mu}})$ is contained in the intersection $$\bigcap_{\substack{s_{i_1},\dots,s_{i_k}\\ v_{i_1},\dots,v_{i_k}}}\widetilde{\mathcal{Y}}_{\substack{s_{i_1},\dots,s_{i_k}\\ v_{i_1},\dots,v_{i_k}}}$$
where $\widetilde{\mathcal{Y}}_{\substack{s_{i_1},\dots,s_{i_k}\\ v_{i_1},\dots,v_{i_k}}}$ is the ideal generated by $$y_{i_1,s_{i_1},v_{i_1}}^{c_{i_1}}y_{i_1+1,s_{i_1+1},v_{i_1+1}}^{c_{i_1+1}}\cdots y_{i_{k},s_{i_k},v_{i_k}}^{c_{i_{k}}}$$
where for each $1\leq j\leq k-1$ we have
$$c_{i_j}+c_{i_j+1}+\cdots+c_{i_{j+1}}= \begin{cases}1,& \text{ if  $s_{i_j}=\cdots=s_{i_{j+1}}=1$}\\
0,& \text{ otherwise. } \end{cases}$$
\begin{proof} This follows from the discussion preceeding the lemma.
\end{proof}
\end{lem}

\begin{example} We consider $r_1=2$ and $r_2=1$ ($i_1=1, i_2=2$) with the two obvious tuples of partitions. In the case $\boldsymbol{\mu}=(\mu_1,\mu_2)$ with $\mu_1=\fontsize{6}{6}\selectfont \Yvcentermath1	\yng(2)$ we get that 
$$(y_{1,2,1}-y_{2,1,1})h_1(\mathbf{y}_{\mu})\in \left \langle y_{1,2,1},y_{2,1,1} \right\rangle$$
and in the other case $\mu_1=\fontsize{6}{6}\selectfont \Yvcentermath1	\yng(1,1)$
$$(y_{1,1,1}-y_{1,1,2})^2h_1(\mathbf{y}_{\mu})\in  \left \langle y_{1,1,1},y_{2,1,1} \right\rangle \cap \left \langle y_{1,1,2},y_{2,1,1} \right\rangle= \left \langle y_{1,1,1}y_{1,1,2},y_{2,1,1} \right\rangle.$$
In the first case we have no restriction on $h_1(\mathbf{y}_{\mu}) $ and in the second case we actually need 
$$h_1(\mathbf{y}_{\mu})\in \left \langle y_{1,1,1}y_{1,1,2},y_{2,1,1} \right\rangle.$$
\end{example}
\begin{cor}\label{maincor}
The dimension of $\varphi_{\boldsymbol{\mu}}(\mathcal{V}_{\boldsymbol{\xi},\gamma,p}\cap \Gamma_{\boldsymbol{\mu}})$ equals the dimension of the space of all homogeneous polynomials $f(\mathbf{y}_{\boldsymbol\mu})$ of degree $-p+|\gamma|-K^{\lambda}_{\boldsymbol{\mu}}$ in the ideal 
$$\bigcap_{\substack{s_{i_1},\dots,s_{i_k}\\ v_{i_1},\dots,v_{i_k}}}\widetilde{\mathcal{Y}}_{\substack{s_{i_1},\dots,s_{i_k}\\ v_{i_1},\dots,v_{i_k}}}$$ symmetric with respect to the exchange of variables $y_{i,r,u}\leftrightarrow y_{i,r,v}$ and satisfying 
$$\text{deg}_{y_{i,r,u}}f(\mathbf{y}_{\boldsymbol\mu})\leq P^{\boldsymbol{\mu},\lambda}_{r,i},\ \forall i,r,u.$$
\begin{proof}
Let $r(\mathbf{y}_{\boldsymbol\mu})$ be a polynomial with $\text{deg}_{y_{i,r,u}}r(\mathbf{y}_{\boldsymbol\mu})\leq P^{\boldsymbol{\mu},\lambda}_{r,i}$, symmetric under the exchange of variables $y_{i,r,u}\leftrightarrow y_{i,r,v}$ and satisfying the property that \eqref{prop3} (with $h_1(\mathbf{y}_{\boldsymbol\mu})$ replaced by $r(\mathbf{y}_{\boldsymbol\mu})$) is contained in \eqref{prop4}. Then, the statement follows if we can show that there exists a function $g(\mathbf{x}_{\gamma})=\frac{f(\mathbf{x}_{\gamma})}{\Delta_{\gamma}}\in \Gamma_{\mu}$ with $f(\mathbf{x}_{\gamma})$ satisfying the properties (1)-(3) from Definition~\ref{mainobj} such that
$$\varphi_{\boldsymbol{\mu}}(g(\mathbf{x}_{\gamma}))=\frac{\Omega_{\boldsymbol\mu}}{\displaystyle \prod_{1\leq i\leq n}\prod_{(r,u)}(y_{i,r,u})^{\lambda_i}}\ r(\mathbf{y}_{\boldsymbol\mu})$$
By the results of \cite{N17a}, we actually obtain the existence of $g(\mathbf{x}_{\gamma}$) with the desired properties except property (3) from Definition~\ref{mainobj}. In the rest of the proof we show that this property can also be fixed. 
Now, we can write $\tilde{f}(\mathbf{x}_{\gamma}):=\prod_{i,r}x_{i,r}^{\lambda_i}f(\mathbf{x}_{\gamma})$ as
$$\tilde{f}(\mathbf{x}_{\gamma})=\tilde{f}_1(\mathbf{x}_{\gamma})+\tilde{f}_2(\mathbf{x}_{\gamma})$$
where $\tilde{f}_1(\mathbf{x}_{\gamma})$ is in the ideal generated by the monomials
\begin{equation}\label{einsideal}(x_{i_1,1}\cdots x_{i_1,r_{i_1}})^{c_{i_1}}(x_{i_1+1,1}\cdots x_{i_1+1,r_{i_1+1}})^{c_{i_1+1}} \cdots (x_{i_k,1}\cdots x_{i_k,r_{i_k}})^{c_{i_k}}\end{equation}
satisfying $c_{i_j}+c_{i_j+1}+\cdots+c_{i_{j+1}}= 1,\ \forall  j\in[1,k)$ and each summand of $\tilde{f}_2(\mathbf{x}_{\gamma})$ is not contained in the aforementioned ideal. This means that for each summand $F$ of $\tilde{f}_2(\mathbf{x}_{\gamma})$ there exists $j\in[1,k)$ and $s_{i_j},\dots,s_{i_{j+1}}$ such that none of the variables
\begin{equation}\label{noneo}x_{i_j,s_{i_j}},x_{i_{j}+1,s_{i_j+1}},\dots ,x_{i_{j+1},s_{i_{j+1}}}\end{equation}
divides $F$.
Since $\varphi_{\boldsymbol{\mu}}(\tilde{f}(\mathbf{x}_{\gamma}))$ and $\varphi_{\boldsymbol{\mu}}(\tilde{f}_1(\mathbf{x}_{\gamma}))$ are both contained in \eqref{prop4}, we have that $\varphi_{\boldsymbol{\mu}}(\tilde{f}_2(\mathbf{x}_{\gamma}))$ is also in \eqref{prop4} 
which forces $\tilde{f}_2(\mathbf{x}_{\gamma})\in\ker \varphi_{\boldsymbol{\mu}}$. If we can show that 
\begin{equation}\label{zzze}\frac{\prod_{i,r}(x_{i,r}^{\lambda_i})^{-1}\tilde{f}_2(\mathbf{x}_{\gamma})}{\Delta_{\gamma}}\in\Gamma_{\mu}\subseteq \mathbb{U}_{\gamma}\end{equation} 
and the nominator of \eqref{zzze} satisfies the properties (1) and (2) from Definition~\ref{mainobj} we obtain that each such polynomial $r(\mathbf{y}_{\boldsymbol\mu})$ is in fact in the image of some element in $\mathcal{V}_{\xi,\gamma,p}\cap \Gamma_{\mu}$. This would finish the proof.
\vspace{0,3cm}

$\bullet \textit{ Symmetry.}$ If we denote by $f^{\sigma}$ the function obtained by applying the symmetric group element $\sigma$ to $f$, we get
$$\tilde{f}_2(\mathbf{x}_{\gamma})-\tilde{f}^{\sigma}_2(\mathbf{x}_{\gamma})$$
is in the ideal generated by the elements \eqref{einsideal}. This is only possible if $\tilde{f}_2(\mathbf{x}_{\gamma})=\tilde{f}^{\sigma}_2(\mathbf{x}_{\gamma})$.

\vspace{0,3cm}
$\bullet \textit{ Specialization.}$ By spacializing $x_{i,1}=x_{i,2}=x_{i\pm 1,1}$ we obtain that
\begin{equation}\label{ssrta}\tilde{f}_1(\mathbf{x}_{\gamma})\Big|_{x_{i,1}=x_{i,2}=x_{i\pm 1,1}=w}=-\tilde{f}_2(\mathbf{x}_{\gamma})\Big|_{x_{i,1}=x_{i,2}=x_{i\pm 1,1}=w}\end{equation}
We consider the case when $i,i\pm1\notin [i_1,i_k]$ or $i,i\pm1\in (i_j,i_{j+1})$ for some $j\in [1,k)$. The remaining cases are proven similarly and we omit the details. Since $\tilde{f}_2(\mathbf{x}_{\gamma})$ is not in the ideal generated by \eqref{einsideal}, there exists $p_1\in [1,k)$ such that 
$$\tilde{f}_2(\mathbf{x}_{\gamma})=R_1(\mathbf{x}_{\gamma})+S_2(\mathbf{x}_{\gamma})$$
where $R_1(\mathbf{x}_{\gamma})$ is independent of $x_{i_{p_1},1},\dots,x_{i_{p_1+1},1}$ and $S_2(\mathbf{x}_{\gamma})$ is divisible by $z$ after specializing $x_{i_{p_1},1}=\cdots=x_{i_{p_1+1},1}=z$. If $i,i\pm1\notin [i_{p_1},i_{p_1+1}]$ we obtain from \eqref{ssrta}  $$R_1(\mathbf{x}_{\gamma})\Big|_{x_{i,1}=x_{i,2}=x_{i\pm 1,1}=w}=0$$
and if $i,i\pm1\in (i_{p_1},i_{p_1+1})$ we obtain that $R_1(\mathbf{x}_{\gamma})$ is divisible by $x_{i,2}$. So collecting all indices $p_1,\dots,p_s\in\{1,\dots,k-1\}$ such that $\tilde{f}_2(\mathbf{x}_{\gamma})$ violates the last condition in Definition~\ref{mainobj} 
we can write by continuing in the above described way
$$\tilde{f}_2(\mathbf{x}_{\gamma})=R_1(\mathbf{x}_{\gamma})+\cdots+R_{s-1}(\mathbf{x}_{\gamma})+x_{i,2} R_s(\mathbf{x}_{\gamma})$$
where $R_j(\mathbf{x}_{\gamma})$ is independent of $x_{i_{p_j},1},\dots,x_{i_{p_j+1},1}$ for $j\leq s$ and is divisible by $z$ after specializing $x_{i_{p_t},1}=\cdots=x_{i_{p_t+1},1}=z$; this property holds for all $t<j$. Moreover, $R_s(\mathbf{x}_{\gamma})$ shows only up if $i,i\pm 1\in (i_{p_s},i_{p_s+1})$ and $$R_j(\mathbf{x}_{\gamma})\Big|_{x_{i,1}=x_{i,2}=x_{i\pm 1,1}}=0,\ \ \forall j<s.$$ Now let $\sigma$ be the element of the symmetric group interchanging the elements $x_{i,1}$ and $x_{i,2}$. On the one hand
$$R_1(\mathbf{x}_{\gamma})-R_1^{\sigma}(\mathbf{x}_{\gamma})$$
is divisible by $z$ after specializing $x_{i_{p_1},1}=\cdots=x_{i_{p_1+1},1}=z$ and on the other hand the term is also independent of $x_{i_{p_1},1},\dots,x_{i_{p_1+1},1}$ since $i\notin [i_{p_1},i_{p_1+1}]$. Hence $R_1(\mathbf{x}_{\gamma})=R_1^{\sigma}(\mathbf{x}_{\gamma})$. Continuing in this way we obtain that $R_j(\mathbf{x}_{\gamma})=R_j^{\sigma}(\mathbf{x}_{\gamma})$ for all $j<s$ and therefore
$$x_{i,2}R_s(\mathbf{x}_{\gamma})=x_{i,1}R_s^{\sigma}(\mathbf{x}_{\gamma}).$$
So $R_s(\mathbf{x}_{\gamma})$ is divisible by $x_{i,1}$ which is only possible if $R_s(\mathbf{x}_{\gamma})=0$ and the claim follows.

\vspace{0,3cm}

$\bullet \textit{ Containement in $\Gamma_{\mu}$.}$ Let $\boldsymbol{\nu}>\boldsymbol{\mu}$ and note that $-\varphi_{\boldsymbol{\nu}}(\tilde{f}_1(\mathbf{x}_{\gamma}))=\varphi_{\boldsymbol{\nu}}(\tilde{f}_2(\mathbf{x}_{\gamma}))$ is contained in the ideal \eqref{prop4} (with the variables corresponding to $\boldsymbol{\nu}$). This implies that $\tilde{f}_2(\mathbf{x}_{\gamma})=\tilde{f}_2(\mathbf{x}_{\gamma})'+\tilde{f}_2(\mathbf{x}_{\gamma})''$ can be written where the first summand belongs to the ideal generated by \eqref{einsideal} and $\tilde{f}_2(\mathbf{x}_{\gamma})''\in \ker \varphi_{\boldsymbol{\nu}}$. But by the choice of $\tilde{f}_2(\mathbf{x}_{\gamma})$ we have $\tilde{f}_2(\mathbf{x}_{\gamma})'=0$ and hence \eqref{zzze} holds.
\vspace{0,3cm}

$\bullet \textit{ First and second property.}$ The first condition of Definitiuon~\ref{mainobj} is obviously satisfied by the nominator of \eqref{zzze} and the second property can be similarly checked. We omit the details.

\end{proof}
\end{cor}
\subsection{} If we denote by $e_d^{(i,r)}(\mathbf{y}_{\boldsymbol\mu})$, $1\leq d\leq m_{i,r}$ the $d$-th elementary symmetric polynomial in the variables $\{y_{i,r,1},\dots, y_{i,r,m_{i,r}}\}$, we can rewrite the statement of Corollary~\ref{maincor} as follows. The dimension of $\varphi_{\boldsymbol\mu}(\mathcal{V}_{\boldsymbol{\xi},\gamma,p}\cap \Gamma_{\boldsymbol\mu})$ equals the dimension of the space spanned by all monomials
$$\prod_{i=1}^n\prod_{r=1}^{r_i}\prod_{d=1}^{m_{i,r}} (e_d^{(i,r)}(\mathbf{y}_{\boldsymbol\mu}))^{C_{d,r,i}}$$
with $$\sum_{d=1}^{m_{i,r}}C_{d,r,i}\leq P^{{\boldsymbol\mu},\lambda}_{r,i},\ \ \sum_{i,r,d}d\cdot C_{d,r,i}=-p+|\gamma|-K^{\lambda}_{\boldsymbol{\mu}},\ \ \sum_{i=i_j}^{i_{j+1}}C_{m_{i,1},1,i}\geq 1, \ \ \forall j\in[1,k).$$
This is exactly the statement of Theorem~\ref{mainthm}.
\begin{rem}\label{partform}
Note that the same proof works for all modules $M_{\boldsymbol{\xi},\lambda}$, $\boldsymbol{\xi}$ normalized,  admitting a presentation as in Lemma~\ref{unminus} where all non-redundant relations in \eqref{gen2} are of the form $(x^-_{i_j,i_{j+1}}\otimes t^{\xi_{i,j}})$, $1\leq i_1<\cdots <i_k\leq n$ satisfying $\lambda_{i_j,i_{j+1}}-\xi_{i_j,i_{j+1}}=1$. 
\end{rem}

\bibliographystyle{plain}
\bibliography{bibfile}

\end{document}